\documentclass{article}

\usepackage[utf8]{inputenc} 
\usepackage[T1]{fontenc}    
\usepackage{amsmath}        
\usepackage{amsfonts}       
\usepackage{amssymb}        
\usepackage{graphicx}       
\usepackage{subcaption}
\usepackage{hyperref}       
\usepackage{amsthm}
\usepackage{xcolor}
\usepackage{easyReview}

\graphicspath{{./Figures/}}

\numberwithin{equation}{section}

\newtheorem{thm}{Theorem}[section]
\newtheorem{prop}[thm]{Proposition}
\newtheorem{cor}[thm]{Corollary}
\newtheorem{lem}[thm]{Lemma}

\theoremstyle{definition}
\newtheorem{defn}[thm]{Definition}
\newtheorem{exa}[thm]{Example}

\theoremstyle{remark}
\newtheorem{rem}[thm]{Remark}

\newcommand{\del}{\partial}

\newcommand{\Conv}{\text{Conv\hspace{1pt}}}
\newcommand{\dist}{\text{dist\hspace{1pt}}}
\newcommand{\orb}{\text{Orb}}

\newcommand{\RR}{\mathbb{R}}
\newcommand{\ZZ}{\mathbb{Z}}
\newcommand{\DD}{\mathbb{D}}
\newcommand{\NN}{\mathbb{N}}

\newcommand{\joliR}{\mathcal{R}}
\newcommand{\joliL}{\mathcal{L}}

\DeclareMathOperator{\length}{Length\vspace{2pt}}

\newcommand{\gr}{\widehat{G_{p,q,r}}}

\title{A Stopping Criterion for the Generation of Conjugacy Classes in Cocompact Triangle Groups}
\author{Jérémy Perazzelli}

\begin{document}
	
	\begin{abstract}
		We give a stopping criterion for the enumeration of all conjugacy classes in cocompact triangle groups up to any geometric length. The enumeration is based on an encoding given by P. Dehornoy and T. Pinsky.
	\end{abstract}
	
	\maketitle
	\tableofcontents

\section{Introduction}\label{sec:intro}

The length spectrum of a hyperbolic surface $X$ is the ordered sequence of lengths $$ \ell(g_1)\leq \ell(g_2)\leq ... $$ of all essential closed geodesics in $X$. We will denote it by $\mathcal{L}(X)$. 

The sequence $\mathcal{L}(X)$ is reminescent of the spectrum of Laplace's operator on hyperbolic surfaces. Indeed, the hyperbolic length of a closed geodesic is analogous to an eigenvalue and a given length may be achieved by multiple geodesics in the same way that eigenvalues may have some multiplicity. Huber's Theorem (\cite{bus} p.229) says that two hyperbolic surfaces $X_1$ and $X_2$ have the same length spectrum if and only if they have the same spectrum for the Laplace operator. The relation between these two spectra is given by Selberg's trace formula.

Triangular surfaces correspond to a certain class of hyperbolic metrics on the thrice punctured sphere. Indeed, for integers $p,q,r$ satisfying $$\frac{1}{p} + \frac{1}{q} + \frac 1 r < 1,$$ the triangle group $G_{p,q,r}$ (defined in Section \ref{sec:intro}) acts on $\DD^2$. Denote by $F$ the set of points with non-trivial stabilizer. Then the quotient space $X_{p,q,r}:=(\DD\setminus F)/G_{p,q,r }$ inherits a hyperbolic structure which is compact when passing through its closure. We will refer to these spaces as cocompact triangular surfaces or only as triangular surfaces. Remark that this is a slight abuse of notation since in most cases, these surfaces are in fact orbifolds. The standard notion of a triangular surface is an hyperbolic surface $\Sigma = \DD/\Gamma$, where $\Gamma\triangleleft G_{p,q,r}$ for some integers $p,q,r$ as above. 

It is a classical fact that non-essential closed geodesics on $X_{p,q,r}$ are in bijection with conjugacy classes of hyperbolic elements in $G_{p,q,r}$. In fact, more is given than a mere bijection: the length of a closed geodesic may be recovered from its corresponding conjugacy class. Hence, the study of the length spectrum of triangular surfaces is completely analoguous to the study of hyperbolic conjugacy classes in triangle groups. To avoid redundancy, in the following, we will refer to hyperbolic conjugacy classes only as conjugacy classes. We also note that there are not parabolic conjugacy class and only three elliptic conjugacy classes given by rotation around every one of the three punctures.

In the article \cite{DP14}, P. Dehornoy and T. Pinsky describe an encoding of all hyperbolic conjugacy classes in any cocompact triangle group. This encoding gives a unique representative for every conjugacy class (but one) in terms of words in the alphabet on the letters $a$ and $b$. Here, we reproduce part of Table 1 from \cite{DP14} dealing with triangular groups for which $p\geq 3$.

\begin{center}\label{Table:first_table}
	\begin{tabular}{c|c|c}
		& $u_L$& $v_R$\\
		\hline\hline
		$p\geq 3$ & &\\
		\hline
		$r$ odd & $((a^{p-1}b)^{\frac{r-3}{2}}a^{p-1}b^2)^\infty$ & $(b^{q-1}a)^{\frac{r-3}{2}} b^{q-1}a^2)^\infty$\\
		\hline
		$r$ even  &$((a^{p-1}b)^{\frac{r-2}{2}}a^{p-2}(ba^{p-1})^{\frac{r-2}{2}}b^2 )^\infty$ & $((b^{q-1}a)^{\frac{r-2}{2}}b^{q-2}(ab^{q-1})^{\frac{r-2}{2}}a^2)^\infty$\\
		\hline
	\end{tabular}
\end{center}

This second table gives the only pair of words which represent the same conjugacy class.
\begin{center}\label{Table:exceptional_words}
	\begin{tabular}{c|c|c}
		& $w_L$& $w_R$\\
		\hline\hline
		$p\geq 3$ & &\\
		\hline
	$r$ odd & $((a^{p-1}b)^{\frac{r-3}{2}}a^{p-1}b)^\infty$ & $(b^{q-1}a)^{\frac{r-3}{2}} b^{q-2}a)^\infty$\\
		\hline
		$r$ even  &$((a^{p-1}b)^{\frac{r-2}{2}}a^{p-2} b (a^{p-1}b)^{\frac{r-4}{2}}a^{p-2}b)^\infty$ & $((b^{q-1}a)^{\frac{r-2}{2}}b^{q-2}a (b^{q-1}a)^{\frac{r-4}{2}}b^{q-2}a)^\infty$\\
		\hline
	\end{tabular}
\end{center}

The encoding in Table \ref{Table:first_table} is done through the notion of a \textit{pair of spectacles} (see Section \ref{sec:spectacles}). This makes the relationship between a closed geodesic an its code somewhat clouded. In particular, if $w$ is the word that represents some closed geodesic $g$, there is not way to tell what is the length of $g$.

In this paper, this is the concern we will be working around when $p\geq 3$. Indeed, after having defined a notion of combinatorial length for words (denoted $L$ in the following), we will establish the following theorem in Section \ref{sec:main_thm}.

\begin{thm}\label{thm:bound_on_length_first_time}
	There exists a constant $c = c_{p,q,r}>0$ which depends only on the group $G_{p,q,r}$ for which we have the inequality $$\length(g) \geq cL(\gamma(g))$$ for any closed geodesic $g$.
\end{thm}

 This theorem has the advantage that the constant $c$ may be obtained explicitely using a computer program. One application is a rigourous criterion to compute the length spectrum of a triangular surface up to a certain bound on the geometric length.
 
 Along with Selberg's Trace Formula, such an algorithm may give rise to some consequences on the multiplicity of eigenvalues of the Laplace operator as was already done in \cite{MFb}.
 
 The paper is divided as follows. In Section \ref{sec:the_graph}, we introduce the graph $\gr$. Paths on this graph will be used to represent closed geodesics. Section \ref{sec:spectacles} is dedicated to explaining the general theory of spectacles developped in \cite{DP14}. Through our reading of this article, we've have had some trouble with justifying some of their claims and so we give other proofs of their results. The content of Sections \ref{sec:some_results about angles in} through \ref{sec:almost_unicity} are designed to fill these holes.
 
 In Section \ref{sec:main_thm}, we give the proof of the Theorem stated above and defined the required concept to enumerate the beginning of the length spectrum without repetition.

\section{The graph $\gr$}\label{sec:the_graph}

In this section, we set up some terminology that will be used in the latter.

\begin{defn}
	A triplet $(p,q,r)\in\NN^3$ is said to be \emph{hyperbolic} if $p\leq q\leq r$ and $$\frac{1}{p} + \frac 1 q + \frac 1 r < 1.$$
\end{defn}

The triplet $(p,q,r)$ defines a unique hyperbolic triangle of internal angles $\frac \pi p,\frac \pi q, \frac \pi r$ coming anticlockwise. Up to isometry, let $A_0,B_0$ and $C_0$ be the vertices of such a triangle in $\DD$. We will denote the full triangle group associated to $(p,q,r)$ by $$\overline{G}_{p,q,r}:=\left\langle r_1,r_2,r_3\right\rangle,$$where each $r_i$ is a reflection in one of the side of the triangle $\Delta(A_0,B_0,C_0)$ and by $G_{p,q,r}\leq \overline{G}_{p,q,r}$ the subgroup of orientation preserving isometries. The graph $\gr$ is the graph whose vertices are the orbits $G_{p,q,r}A_0$ and $G_{p,q,r}B_0
$ (we call such vertices \emph{$A$-type and $B$-type} vertices respectively) and whose edges are the images of the segment $[A_0,B_0]$ under $\overline{G}_{p,q,r}$. The components $\DD\setminus\gr$ will be referred to as \emph{face polygons} and when the context is clear, we will not distinguish between a face polygon $P$ and its closure $\overline{P}$.

Throughout the text, intervals on the boundary $\del\DD$ will oriented clockwise, e.g. for $0 < \theta - \phi< 2\pi$, the interval $(e^{i\theta}, e^{i\phi}]$ represents points $e^{it}$ for $\theta > t \geq \phi$. We will often refer to the cyclic order on $S^1$; by this, we will always mean the \emph{anticlockwise} order, e.g. for a point $\xi \in (e^{i\theta}, e^{i\phi})$, we have $e^{i\phi} < \xi < e^{i\theta}$.

\section{Spectacles}\label{sec:spectacles}

We will describe the tools of \cite{DP14} through which an almost unique encoding of every conjugacy class in $G_{p,q,r}$ is given. This encoding relies on the notion of an \emph{accurate pair of spectacles}. Starting from a fixed base point $P\in \DD$ and an endpoint $\xi\in\del\DD$, a pair of spectacles gives a canonical way to go from $P$ to $\xi$ along $\gr$. For this section, we fix two adjacent vertices $A_0$ and $B_0$ of $\gr$.

\begin{defn}\label{def:base_pair_of_spectacles}
	Denote by $r_{A_0}$ (resp. $r_{B_0}$) the clockwise (resp. anticlockwise) hyperbolic rotation by the angle $\frac{2\pi }{p}$ around $A_0$ (resp. $\frac{2\pi}{q}$ around $B_0$). A \emph{base pair of spectacles} at $A_0,B_0$ is the choice of two intervals $I_{A_0\to B_0} = [\xi_{A_0,B_0}, r_{A_0}(\xi_{A_0,B_0}))$ and $I_{B_0\to A_0} = [\xi_{B_0,A_0}, r_{B_0}^{-1}(\xi_{B_0,A_0}))$, for some points $\xi_{A_0,B_0},\xi_{B_0,A_0}\in\del\DD$. We will often write $S = \{I_{A_0\to B_0}, I_{B_0\to A_0}\}$ for such a pair of spectacles. For every edge $[A,B]$ of $\gr$, a base pair of spectacles yields an \emph{associated} pair of spectacles $gS = \{I_{A\to B}, I_{B\to A}\} $ where $g$ is the unique element of $G_{p,q,r}$ that maps the edge $[A_0,B_0]$ to $[A,B]$, so that $I_{A\to B} = g(I_{A_0\to B_0})$ and $I_{B\to A} = g(I_{B_0\to A_0} )$.
\end{defn}

Suppose that we are given a point $\xi\in\del\DD$, a base pair of spectacles $S$, and some initial vertex $A^0$ of type $A$. We may construct a canonical path $\gamma_{A^0\to \xi}^S$ in the following way.

\begin{enumerate}
	\item Consider all edges $[A^0, B_j]$, where $B_j = r_{A_0}^j(B_0)$ for $ 0 \leq j \leq p-1$. The base pair of spectacles $S$ yields $p$ intervals $I_{A^0\to B_j}$. By construction, these intervals cover $\del\DD$ and are pairwise disjoint. Let $j$ be such that $I_{A^0,B_j}\ni\xi$ and define $B^1 := B_j$. \\
	\item For $n\geq 1$, construct points $A^{2n}$ and $B^{2n+1}$ inductively by considering the vertices $B_k$ of type $B$ adjacent to $A^{2n}$ (or the vertices $A_k$ of type $A$ adjacent to $B^{2n+1}$)
	 and choosing the vertex $B^{2n+1} = B_k$ (or $A^{2n+2} = A_k$) for the only $k$ that satisfies $\xi \in I_{A^{2n}\to B_k}$ (or $\xi \in I_{B^{2n+1}\to A_k}$).
\end{enumerate}

The path $\gamma = \gamma_{A^0\to \xi}^S =  (A^0, B^1, A^2,...)$ is the resulting path. At every vertex $A^{2n}$ for $n\geq 1$, denote by $\theta_{2n}$ the counterclockwise angle between $[A^{2n}, B^{2n-1}]$ and $[A^{2n}, B^{2n+1}]$. By construction, this defines a number $ e_{2n} \equiv \frac{p \theta_{2n}}{2\pi}\mod p$. In an analogous way, one may define $f_{2n+1} \equiv \frac{q \theta_{2n+1}}{2\pi }\mod q$ for the clockwise angle $\theta_{2n+1}$ between $[B^{2n+1}, A^{2n}]$ and $[B^{2n+1}, A^{2n+2}]$ and $n\geq 0$. This yields an encoding of $w(\gamma)= b^{e_0}a^{f_1}b^{e_2}a^{f_2}...$ in the (mono)-infinite monoid generated by powers of $a$ and powers of $b$ where $0\leq e_{2n}\leq p-1$ and $0\leq f_{2n+1}\leq q-1$. Remark however that, in general, the path $\gamma$ may not converge to $\xi$.
\begin{defn}
	Let $P$ be a vertex of $\gr$ and $\xi \in\del\DD$. The path $\gamma_{P\to\xi}^S$ is \emph{$S$-admissible} if it converges to $\xi$. 
\end{defn}

When an $S$-admissible path starting at $P$ exists, it must be unique (c.f. Lemma 3.5 in \cite{DP14}). For two $S$-admissible paths $\gamma_1,\gamma_2$ that start with the same edge, $\gamma_1$ is on the left of $\gamma_2$ after they diverge if and only if $w(\gamma_1)\leq w(\gamma_2)$ with respect to the lexicographic order (c.f. Lemma 3.6 in \cite{DP14}). More precisely, since $\gamma_1$ and $\gamma_2$ are $S$-admissible and start with the same edge, say $[P,P']$, both of their endpoints belong to $I_{P\to P'} = [\xi_1,\xi_2)$. To say that $\gamma_1$ is on the left of $\gamma_2$ is to say that the endpoint $\eta_1$ of $\gamma_1$ and the endpoint $\eta_2$ of $\gamma_2$ satisfy $\xi_2 < \eta_2 < \eta_1 \leq \xi_1$. Remark that by the unicity of $S$-admissible paths, once those two paths diverge, they may never cross again.

The notions defined above also extend to bi-infinite paths. More precisely, given $\eta,\xi\in\del\DD$, a path $\gamma_{\eta\to\xi} = (...,P^{-2},P^{-1},P^0,P^1,P^2,...)$ is $S$-admissible if all its truncations are $S$-admissible paths of the form $\gamma_{P^j\to \xi}$ and $P^j\longrightarrow\eta$ as $j\to-\infty$. This defines an encoding of bi-infinite paths in bi-infinite words.

Among all pair of spectacles, we will want to consider the ones for which every pair of points on the boundary may be joined by an $S$-admissible paths.

\begin{defn}\label{def:biinfiniteadpaths}
	A pair of spectacles $S$ is sub-accurate if for every vertex $P$ and every pair of distinct points $\eta,\xi\in\del\DD$, there exists an $S$-admissible path that joins $P$ to $\xi$ and an $S$-admissible path that joins $\eta$ to $\xi$. 
\end{defn}

This definition is modelled on Definition 3.7 in \cite{DP14}, but our notion of a \textit{sub-accurate} pair of spectacles is their definition of an \textit{accurate} pair of spectacles. We will explain this distinction in Remark \ref{rem:after_strong_accuracy_definition}. In that paper, a concrete pair of spectacle, called $S^f$ (see Section \ref{sec:explicit_pair_of_spec}) is defined and it is shown that it is accurate, i.e. the authors show (c.f. Lemma 3.13) that every path $\gamma^{S^f}_{P\to\xi}$ as above converges and construct (c.f. Lemma 3.14) an $S^f$-admissible path from $\eta$ to $\xi$ whenever $\eta,\xi\in\del\DD$ are distinct. A couple of things should be said about those joinings. 

We first point out that the proof of Lemma 3.14 doesn't rely on the definition of $S^f$; the admissible path they construct is obtained by a limiting process involving mono-infinite admissible paths $\gamma^{S^f}_{Q\to \xi}$. This means that in Definition \ref{def:biinfiniteadpaths}, we may remove the condition that $S$-admissible paths joining arbitrary points on the boundary exist, since this follows from the other condition.

Now, let $P = A^0$ be a vertex of $\gr$ (which is taken to be of type $A$ for simplicity) and $\xi\in\del\DD$, denote by $B_{\xi}$ the Busemann function based at $\xi$ relative to the point $A^0$ and $\gamma^{S^f}_{A^0\to\xi} = (A^0, B^1, A^2,...)$. For the proof of Lemma 3.13, the author claim that $\mathcal{B}_{\xi}(B^{{2j-1}}) - \mathcal{B}_\xi( B^{2j+1})>0$ for all $j\geq 0$ and refer to a figure. The goal of Section \ref{sec:monotony_of_bus_fun} is to give a proof of this statement. 




We now discuss how $S$-admissible paths may be characterized more easily by using the lexicographical order.

\begin{defn}
	Let $S$ be a base pair spectacles. Define $S' = \{I_{A_0\to B_0}' = (\xi_{A_0,B_0}, r_{A_0}(\xi_{A_0,B_0})],I_{B_0\to A_0}' = (\xi_{B_0,A_0}, r_{B_0}(\xi_{B_0,A_0})]\}$. We call $S'$ the \emph{dual base pair of spectacles to $S$}. As in Definition \ref{def:base_pair_of_spectacles}, we get intervals \emph{dual associated spectacles} $I'_{A\to B}$ and $I'_{B\to A}$ for every edge $[A,B]$ in $\gr$. The intervals in a regular pair of spectacles were closed on the left and open on the right. In a dual pair of spectacles, intervals are open on the left and closed on the right.
\end{defn}

Dual spectacles were not defined in \cite{DP14}, but we feel that they shed light on the following definition. As regular pair of spectacles, dual pair of spectacles define canonical paths $\gamma^{S'}_{P\to \xi}$ joining a vertex $P$ to $\xi\in\del\DD$. When such paths converge, we will say that they are $S$-dually-admissible or $S'$-admissible.

\begin{defn}
	A pair of spectacles $S$ is said to be \emph{accurate} if it is sub-accurate and the paths $\gamma_{A_0\to r_{A_0}(\xi_{A_0,B_0})}^{S'}$ and $\gamma_{B_0\to r_{B_0}(\xi_{B_0,A_0})}^{S'}$ are $S$-dually-admissible. 
\end{defn}

\begin{rem}\label{rem:after_strong_accuracy_definition}
	As stated above, our notion of accuracy is stricly stronger than the one in \cite{DP14}. We ask that accurate pairs of spectacles have their rightmost path to be dually admissible. Since $S$-admissibility is hard to check in practice, we will resort to a combinatorial criterion to check if a given word is admissible. This criterion will make use of a leftmost path and a rightmost path. However, in \cite{DP14}, this rightmost path is ill-defined and its existence is absolutely necessary for this criterion to work. This justifies our stronger definition for accuracy.
\end{rem}

\begin{defn}\label{def: extreme-most_words}
	Let $S = \{I_{A_0\to B_0} = [\xi_{A_0,B_0}, r_{A_0}(\xi_{A_0,B_0})),I_{B_0\to A_0} = [\xi_{B_0,A_0}, r_{B_0}^{-1}(\xi_{B_0,A_0}))\}$ be an accurate pair of spectacles. We define $v_L := w\left(\gamma_{A_0\to \xi_{A_0,B_0}}^{S}\right)$, $v_R := w\left(\gamma_{A_0\to r_{A_0}(\xi_{A_0,B_0})}^{S'}\right)$, $u_L := w\left(\gamma_{B_0\to\xi_{B_0,A_0}}^S\right)$ and $u_R := w\left(\gamma_{B_0\to r_{B_0}^{-1}(\xi_{B_0,A_0})}^{S'}\right)$.
\end{defn}

Here, we point out that our definitions for $(v_L,v_R)$ and $(u_L,u_R)$ are interchanged with respect to the ones in \cite{DP14}. This is justified by what seems to be a typographical error in \cite{DP14}. Indeed, the words $u_L$ and $u_R$ correspond to paths starting with a point of type $A$, so they should start with a letter $b^f$ which is not the case (see Table 1 in \cite{DP14}). Besides, our definition of the right-most paths $v_R$ and $u_R$ are fundamentally different than in \cite{DP14}. 

Using the argument in Lemma 3.6 of \cite{DP14}, we can show that a bi-infinite path $\gamma = (...,P^{k-1}, P^k, ...)$ is $S$-admissible if and only if for all $n\in\ZZ$, the truncation $\gamma_n = (P^n, P^{n+1},...)$ satisfy $s_L\leq_{lex} w(\gamma_n) <_{lex} s_R$ for $s = v$ when $P^n$ is of type $A$ and $s = u$ when $P^n$ is of type $B$. This is essentially the content of Lemma 3.10 in this article. Indeed, their definition asks for $u_R$ to be the $S$-admissible path joining $A_0$ to $r_{A_0}(\xi_{A_0,B_0})$. However, this path is obtained from a rotation of $u_L$, so they both have the same code. These complete the concerns that were raised in Remark \ref{rem:after_strong_accuracy_definition}. We resume this discussion and the part that is most relevant for us in the following proposition.

\begin{prop}\label{prop: admissible_iff_some_inequalities}
	Let $S$ be an accurate pair of spectacles. Then a bi-infinite path $\gamma = (...,P^{k-1}, P^k, ...)$ is $S$-admissible if and only if for all $n\in\ZZ$, the truncation $\gamma_n = (P^n, P^{n+1},...)$ satisfy $s_L\leq_{lex} w(\gamma_n) <_{lex} s_R$ for $s = v$ when $P^n$ is of type $A$ and $s = u$ when $P^n$ is of type $B$.
\end{prop}


\section{Some results about angles in $\gr$}\label{sec:some_results about angles in}

From this point on, our goal is to describe the said pair of spectacles $S^f$ defined in \cite{DP14} and to show that it is accurate. The pair of spectacles $S^f$ will be defined in Section \ref{sec:explicit_pair_of_spec} and its accuracy will be established in Section \ref{sec:monotony_of_bus_fun}. This section is a collection of technical lemmas that will be used in further sections and should maybe be read as references to those lemmas are made.

\begin{figure}
	\centering
	\includegraphics{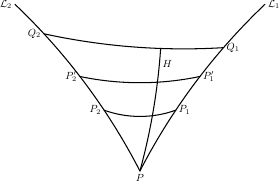} 
	\caption{The three triangles $\Delta_1$, $\Delta_2$, $\Delta_3$ defined in Lemma \ref{lem:conditionangle}.}
	\label{fig:triangle}
\end{figure}

Let $\joliL_1$ and $\joliL_2$ be two hyperbolic lines that meet with angle $0 < \alpha < \pi$ in a point $P$ and consider, for $\psi_1>0$, the isosceles triangle $\Delta_1$ emanating from $P$ with interior angles $\alpha,\psi_1,\psi_1$ (if it exists). We denote by $P_i$ the vertex of $\Delta_1$ on $\joliL_i$ and define $\ell:=\length([P,P_1])$. Then, for $i=1,2$, consider a point $Q_i$ lying of $\joliL_i$ further than $P_i$. Construct the triangle $\Delta_3 = \Delta(P,Q_1,Q_2)$. Denote by $\psi_1'$ and $\psi_2'$, respectively, the angles $\angle(P,Q_1,Q_2)$ and $\angle(P,Q_2,Q_1)$ (see Figure \ref{fig:triangle} for reference). Remark that by the Gauss-Bonnet Theorem, we must have $\psi_1'\leq\psi_1$ or $\psi_2'\leq\psi_1$, but that one of these two inequalities may be false. Finally, let $H$ be the height of $\Delta_3$ based at $P$ and $M$ be its other endpoint, see Figure \ref{fig:triangle}.

\begin{lem}\label{lem:conditionangle}
	With the notation above, and assuming that $\psi_2'<\psi_1'$ (this ensures, in particular, that $\beta < \frac{\alpha}{2}$), let $h := \length(H)$ and $\beta$ be the angle made by $\joliL_1$ and $H$. Define $s_i :=\length[P,Q_i]$. Then $\psi_1'<\psi_1$ if and only if
	\begin{equation}\label{eq:lemtriangles}
		\cosh s_1\bigg(\frac{\tanh s_1}{\tanh s_2} - \cos\alpha\bigg)> (1-\cos\alpha)\cosh \ell.
	\end{equation}
\end{lem}

\begin{proof}
	Consider the isosceles triangle $\Delta_2$ having vertices $P_1'$ on $\joliL_1$ and $P_2'$ on $\joliL_2$ such that its internal angles are $\alpha,\psi_1',\psi_1'$. Then $\psi_1' < \psi_1$ if and only if $\Delta_1\subsetneq\Delta_2$. If we let $\ell':=\length[P,P_1']$, then this is equivalent to $\cosh\ell < \cosh\ell'$. All trigonometric formulas can be found in \cite{bus} on page 454. Using the second law of cosines, we have that 
	\begin{equation}\label{eq:cosh}
			\cosh\ell' = \frac{\cos\psi_1' + \cos\alpha\cos\psi_1'}{\sin\psi_1'\sin\alpha} = \cot\psi_1'\bigg(\frac{1+\cos\alpha}{\sin\alpha}\bigg).
	\end{equation}
	
	Since $\Delta(P,M,Q_1)$ is a right angled triangle, $\cosh s_1 = \cot\beta \cot\psi_1'$. The RHS of \eqref{eq:cosh} becomes
	\begin{equation}\label{eq:tanbeta1}
		\tan\beta\cosh s_1\bigg(\frac{1+\cos\alpha}{\sin\alpha}\bigg).
	\end{equation}
	
	Since $\Delta(P,Q_2,M)$ also is a right-angled triangle, we have the two following equations
	\begin{equation}\label{eq:s1}
		\cos(\alpha-\beta) = \tanh h \coth s_1
	\end{equation}
	\begin{equation}\label{eq:s2}
		\cos\beta = \tanh h\coth s_2.
	\end{equation}
	
	From equations \eqref{eq:s1} and $\eqref{eq:s2}$, we have 
	\begin{equation}\label{eq:tanbeta2}
		\frac{\tanh s_1}{\tanh s_2} = \cos\alpha + \sin\alpha \tan\beta.
	\end{equation}
	Combining equations \eqref{eq:cosh}, \eqref{eq:tanbeta1} and \eqref{eq:tanbeta2}, the inequality $\cosh\ell < \cosh \ell ' $ is seen to be equivalent to 
	\begin{align*}
		\cosh\ell  &< \frac{1}{\sin \alpha}\left( \frac{\tanh s_1}{\tanh s_2}-\cos\alpha\right)\left(\frac{1+\cos\alpha}{\sin\alpha}\right)\cosh s_1\\
		&=\left(\frac{\tanh s_1}{\tanh s_2} - \cos \alpha \right ) \frac{\cosh s_1}{1-\cos\alpha},
	\end{align*}
	which yields the desired inequality.
\end{proof}

\begin{figure}
	\centering
	\includegraphics{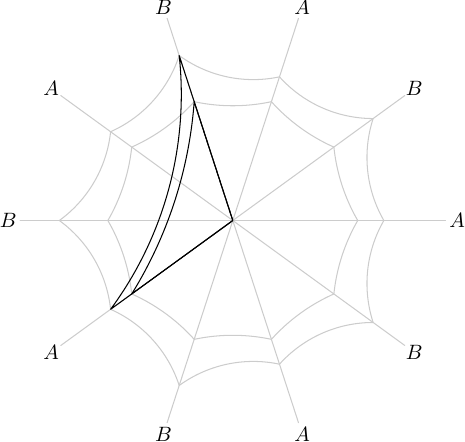} 
	\caption{The two triangles in bold correspond to $\Delta_1$ and $\Delta_3$. The big gray polygon corresponds to a case where $p_0 < q$ and $r=5$. The smaller one corresponds to $q = p_0$ and $r=5$.
		\label{fig:triangle_in_pol}}
\end{figure}

Figure \ref{fig:triangle_in_pol} shows how the previous situation occurs in the study of $\gr$.

First suppose that $3\leq p_0 \leq r_0$.  Let $P\subseteq \DD$ be a face polygon with center $C$ and fix one of its edges $e=[x,y]$ oriented conveniently with the standard orientation on $P$. Starting with $e$, label all edges from $1$ to $2r$ counterclockwise. Starting from $y$, label the vertices counterclockwise from $1^+$ to $2r^+$. All indices are taken (mod $2r$). Now, let $q\in\NN_{\geq 3}$ be such that $(p_0,q,r_0)$ is a hyperbolic triplet. Define $\nu(q):=\nu_q:=\angle (2r)^+, (2r-1)^+, (r+1)^+$, $\mu(q):=\mu_q:=\angle r^+, (r+1)^+, (2r-1)^+$, $\phi(q):=\phi_q:=\angle (2r-1)^+, (r+1)^+, C$ and $\psi(q):=\psi_q:=\angle (r+1)^+,(2r-1)^+, C$. Here, there are four possible configurations depending on if $e=[A,B]$ or $[B,A]$ and on the parity of $r$:
\begin{enumerate}
	\item if $e = [A,B]$ and $r$ is odd, then $(2r-1)^+$ if of type $B$ and $(r+1)^+$ is of type $A$;\label{case:1}\\
	\item if $e=[A,B]$ and $r$ is even, then $(2r-1)^+$ and $(r+1)^+$ are of type $A$;\label{case:2}\\
	\item if $e=[B,A]$ and $r$ is odd, then $(2r-1)^+$ is of type $A$ and $(r+1)^+$ is of type $B$;\label{case:3}\\
	\item if $e=[B,A]$ and $r$ is even, then $(2r-1)^+$ and $(r+1)^+$ are of type $B$.\label{case:4}
\end{enumerate}

We define $\alpha:=\angle (2r-1)^+, C, (r+1)^+$. In any of the four preceding cases, we have $\alpha = \pi - \frac{2\pi}{r}$.

\begin{figure}[h]
	\centering
	\begin{subfigure}[b]{0.45\textwidth}
		\includegraphics[width=\textwidth]{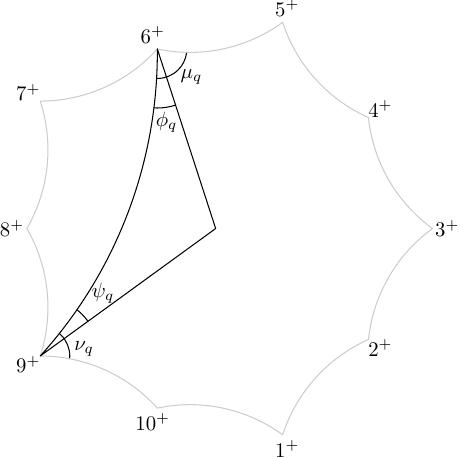}
		\caption{Case \ref{case:1}}
		\label{fig:case1}
	\end{subfigure}
	\hfill
	\begin{subfigure}[b]{0.45\textwidth}
		\includegraphics[width=\textwidth]{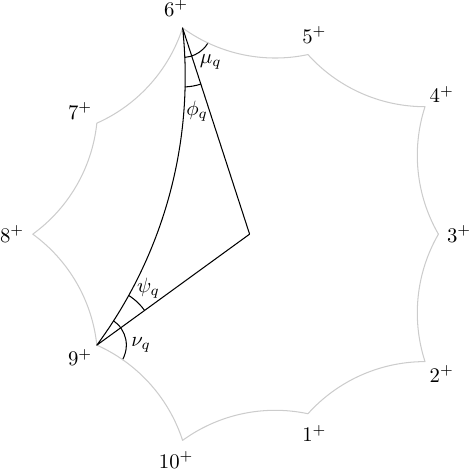}
		\caption{Case \ref{case:3}}
		\label{fig:case3}
	\end{subfigure}
	\caption{Cases \ref{case:1} and \ref{case:3} can be obtained one another by hyperbolic isometries.}
	\label{fig:main}
\end{figure}

\begin{rem}\label{rem:symmetry_of_cases}
	As shown in Figures \ref{fig:case1} and \ref{fig:case3}, cases \ref{case:1} and \ref{case:3} are in fact the same. Indeed, if one reflects case \ref{case:3} in the axis of $[(2r-1)^+,(r-2)^+]$ and rotates it clockwise by $\pi-\frac {2\pi}{r}$ around $C$, we get the picture of case \ref{case:1} where $\phi_q$ plays the role of $\psi_q$ and the $\mu_q$ plays the role of $\nu_q$ and \textit{vice versa}.
\end{rem}

\begin{lem}\label{lem:angle_decreases_odd_case}
	Let $3\leq p_0\leq r_0$ be as above. In case \ref{case:1}, $\psi_q\leq \phi_q\leq\phi_{p_0}$. In case \ref{case:3}, $\phi_q\leq \psi_q \leq \psi_{p_0}$. 
\end{lem}

\begin{proof}
	
	First note that by Remark \ref{rem:symmetry_of_cases}, the two things we need to show are in fact the same. Therefore, we will only prove the first statement.
	
	If $q=p_0$, we have $\phi_q = \phi_{p_0}$, so we will suppose $q> p_0$. If we triangulate $P$ by joining its vertices to $C$, we get isometric images of the fundamental triangle with internal angles $\frac{\pi}{p},\frac \pi q$ and $\frac \pi r$. The side $[C,(r+1)^+]$ is opposite to an angle $\frac \pi q$ and side $[C,(2r-1)^+]$ is opposite to an angle $\frac \pi p$. Hence, we will write $b_q:=\length([C,(r+1)^+])$ and $a_q = \length([C,(2r-1)^+])$. Note that, by the sine law, $a_q > b_q$ so that $\psi_q < \phi_q$. In reference to Lemma \ref{lem:conditionangle}, $\psi_q$ plays the role of $\psi_2'$ and $\phi_q$ plays the role oh $\psi_1'$. By the aforementioned Lemma, we need to show that
	\begin{equation}\label{eq:amontrerlemanglepoly}
		\cosh b_q \bigg(\frac{\tanh b_q}{\tanh a_q} - \cos\alpha\bigg)>(1-\cos\alpha)\cosh b_p .
	\end{equation}
	
	In the following, to make the notation less cumbersome, we will write $p$ for $p_0$ and $r$ for $r_0$. Using the fundamental triangles and by the sine and cosine laws, we have
	\begin{equation}\label{eqthth}
		\frac{\tanh b_q}{\tanh a_q} = \frac{\cos\tfrac{\pi}{p} +\cos\tfrac{\pi}{q}\cos\tfrac \pi r}{\cos\tfrac \pi q + \cos\tfrac \pi p \cos \tfrac \pi r}
	\end{equation}
	and
	\begin{equation}\label{eq:coshbp}
		\frac {\cosh b_p} {\cosh b_q} = \frac{\cos\tfrac \pi p + \cos \tfrac \pi p \cos \tfrac \pi r}{\cos\tfrac \pi q + \cos \tfrac \pi p \cos \tfrac \pi r}.
	\end{equation}
	
	Hence \eqref{eq:amontrerlemanglepoly} becomes
	\begin{equation}\label{eq:amontrerangle}
		\cos\tfrac \pi p + \cos \tfrac \pi q \cos \tfrac \pi r - (\cos\tfrac \pi q + \cos \tfrac \pi p \cos \tfrac \pi r)\cos\alpha > (1-\cos\alpha)(\cos\tfrac \pi p + \cos \tfrac \pi p \cos \tfrac \pi r).
	\end{equation}
	Now, computing the derivative of the LHS of the latter with respect to $\cos\tfrac{\pi}{q}$ yields the term $\cos\tfrac \pi r - \cos\alpha$ which is positive since $\alpha = \pi - \tfrac{2\pi}{r}$ and $r\geq 4$. Hence the left hand side of \eqref{eq:amontrerangle} is a strictly increasing function of $q$ and the inequality becomes an equality precisely when $q = p$.
\end{proof}

\begin{lem}\label{lem:generic_case_smaller_than_pi2}
	Let $(p_0,q,r_0)$ be as above, but suppose that $(p_0,q,r_0)\not = (3,3,4)$ or $(3,4,4)$. Then, in all cases, $\mu_q < \frac{\pi}{2}$ and $\nu_q < \frac{\pi}{2}$.
\end{lem}

\begin{proof}
	As we have explained above, because of the symmetries, we only need to check the proposition for case \ref{case:1}, \ref{case:2} and \ref{case:4}. We will begin with the case \ref{case:1}.
	
	First remark that $\mu_q =\phi_q+\frac{\pi}{p}$ and $\nu_q = \psi_q + \frac \pi q$. Since, $b_q\leq a_q$, by the sine law, we have that $\psi_q\leq\phi_q$. Hence, $\nu_q = \psi_q +\frac\pi q\leq \phi_q + \frac\pi p=\mu_q$. So we have reduced ourselves to show that $\mu_q < \frac\pi 2$.
	
	By Lemma \ref{lem:angle_decreases_odd_case}, we have that $\phi_q\leq\phi_p$ and the angle $\phi_p$ corresponds to the more symmetric case $(p,p,r)$. In this case, the triangle $\Delta((2r-1)^+,C, (r+1)^+)$ is isosceles (see the smaller triangle in bold in Figure \ref{fig:triangle_in_pol}). If we cut along its height, we may obtain a right angled triangle of internal angles $\frac
	\alpha 2$ and $\phi_p$ and hypotenuse of length $b_p$. Formulas from hyperbolic trigonometry give that
	\begin{equation}\label{eq:hyp_tri_right_angle1}
		\cosh b_p = \cot\frac \alpha 2 \cot \phi_p.
	\end{equation}
	Now, recall that $\cosh b_p = \cot\frac{\pi}{p}\left(\frac{1+\cos\frac{\pi}{r}}{\sin\frac \pi r}\right)$ and $\frac \alpha 2 = \frac \pi 2 - \frac \pi r$ so that \eqref{eq:hyp_tri_right_angle1} is equivalent to 
	\begin{equation}\label{eq:lower_bound1}
		\cot \phi_p = \cot\frac \pi p \left(\frac{\cos \frac \pi r}{1-\cos \frac \pi r}\right).
	\end{equation}
	Remark that \eqref{eq:lower_bound1} is an increasing function of $p$ and $q$ as a product of two increasing functions. 
	
	Since we are considering case \ref{case:1}, $r$ is odd. So, for $(p_0,q,r_0)$ to be a hyperbolic triplet, we must have $r\geq 5$. Equation \eqref{eq:lower_bound1} then satisfies
	
	\begin{align*}
		\cot\phi_p =  \cot\frac \pi p \left(\frac{\cos \frac \pi r}{1-\cos \frac \pi r}\right) \geq \cot \frac \pi 3\left(\frac{\cos \frac\pi 5}{1 - \cos \frac \pi 5}\right)=:L.
	\end{align*}
	
	Taking $\cot^{-1}$ yields, and considering $\mu_q$, we obtain
	\begin{align*}
		\mu_q = \phi_q + \frac \pi p\leq \phi_p + \frac \pi 3 \leq \cot^{-1}L+\frac \pi 3 \approx 1.435... < \frac \pi 2.
	\end{align*}

	For cases \ref{case:2} and \ref{case:4}, the triangle $\Delta((2r-1)^+,C,(r+1)^+)$ is always isosceles. Hence $\phi_q \psi_q$ and $\mu_q = \nu_q$. By the law of sines, to decrease $q$ amounts to decrease the length of the sides of equal length of this triangle. This increases the angles $\phi_q $. Hence, $\phi_q \leq \phi_p$. 
	
	As in case \ref{case:1}, we obtain a right angled triangle by cutting $\Delta((2r-1)^+,C,(r+1)^+)$ in two along its height. The same formulas apply and we end up with 
	
	\begin{align*}
		\cot\phi_p =  \cot\frac \pi p \left(\frac{\cos \frac \pi r}{1-\cos \frac \pi r}\right).
	\end{align*}
	Here, since we don't consider the triplets $(3,3,4)$ and $(3,4,4)$ which are the only two triplets with $r = 4$, we may suppose as in case \ref{case:1} that $r\geq 5$ and we can conclude in the same way.
\end{proof}



\section{Order relations on some paths}

In this section, we will introduce some order relations used on some paths in $\DD$. In the following, we will ofter consider paths $\alpha$ on $\gr$ as a sequence $(\alpha_j)_{j\in\ZZ}$, where the $\alpha_j$'s are the vertices of $\alpha$ and geodesic segments $[\alpha_j,\alpha_{j+1}]$ are its (oriented) edges. Note that this ``parametrization'' is not unique. If $\alpha$ tends to a point $\eta$ in the $\pm\infty$ direction, we will write $\alpha_{\pm\infty} = \eta$.

\begin{defn}
	Let $\eta$ and $\xi\in\del\DD$. A path $\alpha$ joining $\eta$ to $\xi$ is a \emph{embedded proper path} if it is a topological embedding $\alpha:\RR\to \DD$ that satisfies $\alpha_{-\infty}=\eta$ and $\eta_\infty = \xi$. The set of such paths will be denoted by $\Sigma_{\eta,\xi}$.
\end{defn}

Two emdebbed proper paths will be considered equivalent if there exists an increasing reparametrization $s:\RR\to\RR$ such that $\beta = \alpha\circ s$. Along the text, we will denote by $\mathcal{P}(\gr)$ the set of paths on the graph $\gr$.


The fact that embedded proper paths are topological embeddings enables us to put some structure on $\Sigma_{\eta,\xi}$. Indeed, by the Jordan curve Theorem, the set $\DD\setminus \alpha$ is the disjoint union of two cells and, because of the orientation on $\alpha$, one is naturally on the left, and the other, on the right of $\alpha$. We will denote these cells respectively by $L_\alpha$ and $R_\alpha$. 

\begin{defn}
	If $\alpha,\beta\in\Sigma_{\eta,\xi}$, we define $\alpha\leq \beta$ if $\beta\subseteq \overline{R_\alpha}$, $\alpha\ll \beta$ if $\alpha\leq \beta$ and their sets of edges are disjoint and $\alpha\prec \beta$ if $\beta\subseteq R_\alpha$.
\end{defn}

It is readily seen that $\leq$ is symmetric. For reflexivity, suppose that $\alpha\leq \beta$ and that $\beta \leq \alpha$. Since the statement is purely topological, we may assume that we are working in the square $S = (-1,1)\times(-1,1)$, that $\alpha(t) = (0,t)$ for $t\in(-1,1)$ and that $L_\alpha = (-1,0)\times(-1,1)$ and $R_\alpha = (0,1)\times(-1,1)$.

Our first claim is that $\alpha \leq \beta \implies \alpha\subseteq\overline{L_\beta}$. Recall that $S = L_\beta\sqcup \beta \sqcup R_\beta$. Now, consider, for $n\geq 1$ the rectangle $R^n = (\frac{-1}{n},1)\times (-1,1)$, so that $\beta\subsetneq R^n$ for all $n$. By the Jordan curve Theorem, $\beta$ divides each $R^n$ in $L_\beta ^n\sqcup \beta\sqcup R_\beta^n$, where $L_\beta^n = L_\beta\cap L^n$ and $R_\beta^n = R_\beta\cap R^n = R_\beta$. By construction, we also have that $\bigcap_n R^n = \overline{R_\alpha}$ and $R_\beta = R_\beta^n\subseteq R^n$. Taking the intersection yields $R_\beta\subseteq\overline{R_\alpha}$ which, in its turn, implies that $\overline{R_\beta}\subseteq\overline{R_\alpha}$. Looking first at the complement and then at the closure, we have that $\overline{L_\beta}\supseteq\overline{L_\alpha}\supseteq\alpha$.

Since we also supposed that $\beta \leq \alpha$, we can conclude that $\alpha \subseteq \overline{L_\beta}\cap\overline{R_\beta} = \beta$. A similar argument shows that $\beta\subseteq\alpha$. Without loss of generality, we may assume that $\beta$ is parametrized by $(-1,1)$ in a way that $\lim_{t\to -1^+} \beta(t) = (0,-1)$ and $\lim_{t\to 1^-} \beta(t) =(0,1)$. Since $\alpha = \beta$, we have that for every $t\in(-1,1)$, there exists $s(t)\in(-1,1)$ such that $\alpha(t) = \beta(s(t))$ because $\beta$ is homeomorphism onto its image, $s$ is continous and can be show to be invertible, hence monotonous. The fact that $\alpha(t),\beta(t)\longrightarrow(-1,0)$ as $t\to -1^+$ and $\alpha(t),\beta(t)\longrightarrow(0,1)$ as $t\to 1^-$  implies that $s$ in increasing, so that $\alpha\sim \beta$. This establishes the reflexivity of $\leq$. 

In the same order of idea, $\leq$ is transitive, thus a partial order on $\Sigma_{\eta,\xi}$. By a similar argument, $\ll$ and $\prec$ are strict order relations. The three relations are linked by the following transitivity property:

\begin{itemize}
	\item if $\alpha\leq \beta$ and $\beta\ll \gamma$ or if $\alpha\ll\beta$ and $\beta\leq \gamma$, then $\alpha\ll\gamma$;\\
	\item  if $\alpha\leq \beta$ and $\beta\prec \gamma$ or if $\alpha\prec\beta$ and $\beta\leq \gamma$, then $\alpha\prec\gamma$;\\
	\item if $\alpha\ll \beta$ and $\beta\prec \gamma$ or if $\alpha\prec\beta$ and $\beta\ll \gamma$, then $\alpha\prec\gamma$.
\end{itemize}

\begin{defn}\label{def:landR-convex}
	If $\alpha\in \Sigma_{\eta,\xi}$, we say that $\alpha$ is \emph{left-convex} (resp. \emph{right-convex}) if $\overline{L_\alpha}$ (resp. $\overline{R_\alpha})$ is convex.
\end{defn}

\section{Explicit pair of spectacles}\label{sec:explicit_pair_of_spec}
Let $(p,q,r)$ be a hyperbolic triplet and $(P_j)_{j\in \ZZ}$ be a sequence of face polygons that are adjacent in opposite sides, i.e. for every $j\in \ZZ$, the edges $P_{j-1}\cap P_j$ and $P_j\cap P_{j+1}$ are opposite one another in $P_j$. Let $e_j$ be the edge $P_j\cap P_{j+1}$. To each edge $e_j$, we give the orientation induced by $P_{j+1}$. Accordingly to this orientation, we write $e_j^-$ and $e_j^+$ for the endpoints of $e_j$. 

\begin{lem}\label{lem: infinite_edge_seq}
	With the notation above, define $\omega_L := \bigcup_{j\in\ZZ}[e_j^-,e_{j+1}^-]$ and $\omega_R:=\bigcup_{j\in \ZZ}[e_{j}^+, e_{j+1}^+]$. Then $\omega_L$ is right-convex, $\omega_R$ is left-convex and both path converge to the same limit points. In particular, there is a unique geodesic $g$ contained in $\Conv (\omega_L)\cap \Conv(\omega_R)$ joining the limit points. We also have the inequality $$\omega_L\ll g \ll \omega_R.$$
\end{lem}

\begin{proof}
	To prove that $\omega_L$ is right-convex, we need to prove that for all $j\in\ZZ$, the angle $\theta^j := \angle e_{j+1}^-,e_j^-, e_{j-1}^-$ is $\leq \pi$. We consider $\theta^j = \mu^j +\nu^j$, where $\mu^j := \angle e_{j-1}^-,e_j^-,e_j^+$ and $\nu^j := \angle e_j^+, e_j^-, e_{j+1}^-$ (see Figure \ref{fig:angles_theta_with_edges_e}).
	
	\begin{figure}
		\centering
		\includegraphics{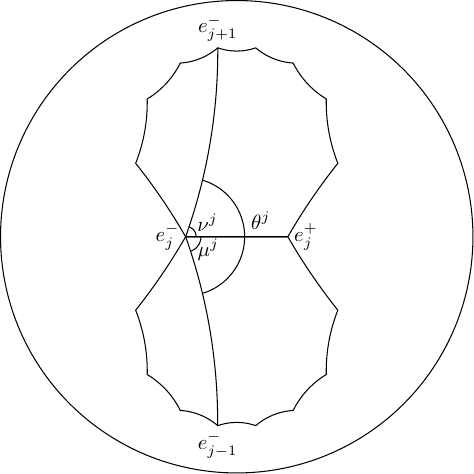} 
		\caption{The angles $\theta^j$ are convex.
			\label{fig:angles_theta_with_edges_e}}
	\end{figure}

	If $(p,q,r)\not = (3,3,4)$ or $(3,4,4)$, we can compare angles $\mu^j$ and $\nu^j$ with the ones defined in Lemma \ref{lem:generic_case_smaller_than_pi2}. Since those are always (i.e. in any of the four cases described) strictly smaller than $\frac \pi 2$, we have $\theta^j< \frac \pi 2+\frac \pi 2 = \pi$.
	
	In the case $(p,q,r) = (3,3,4)$, direct computations show that whether $e_j^-$ is an $A$-type or $B$-type vertex, we always have $\theta^j\approx 2.672... < \pi$.
	
	In the case $(p,q,r) = (3,4,4)$, we have $\theta^j\approx 2.579...$ if $e_j^-$ is of type $A$ and $\theta^j\approx 1.965...$ if it is a $B$-type vertex.
	
	For $\omega_R$, the statement follows by the same reasoning. 
	
	Since $\omega_L$ and $\omega_R$ stay at bounded distance one from another they must converge to the same point at infinity. Since $\omega_L$ is right-convex and $\omega_R$ is left-convex, there must exist a geodesic $g$ satisfying $\omega_L \leq g \leq \omega_R$. Since $\forall j \in \ZZ$, $\theta^j < \pi$, no edge of $\omega_L$ or $\omega_R$ can intersect fully with $g$. Hence, we have $\omega_L \ll g \ll \omega_R$.
\end{proof}





We now review the construction of accurate spectacles done in \cite{DP14}.

\begin{figure}
	\centering
	\includegraphics{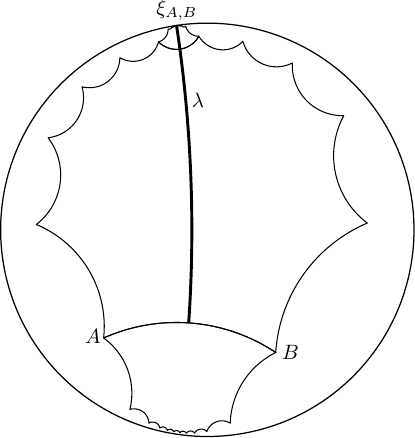} 
	\caption{The half geodesic $\lambda$ plays the role of the spectacles defining geodesic of $I_{A\to B}$ and $I_{B\to A}$.
	\label{fig:spec_defining)geod}}
\end{figure}

Given an oriented edge $[A,B]$, let $P_1$ be the face polygon on its left and $P_0$ be the polygon on its right. Construct a sequence of opposite face polygons $(P_j)_{j\in\ZZ}$ as it was done in the discussion preceding Lemma \ref{lem: infinite_edge_seq}. Then by the said Lemma, there is a unique geodesic $g$ contained in the union $\bigcup_{j\in \ZZ} P_j$. We define $\xi_{A,B}:=\lim_{t\to\infty} g(t)$. We will denote by $\lambda$ the half geodesic joining the intersection point of $[A,B]$ with $g$ to $\xi_{A,B}$ (see Figure \ref{fig:spec_defining)geod}). If the edge is oriented in the opposite direction, i.e. we consider an edge $[B,A]$, then we take $P_1$ as the face polygon on the \emph{right} of $[B,A]$ and $P_0$ as the face polygon on its left. Invoking Lemma \ref{lem: infinite_edge_seq} once again yields a analogous half-geodesic $\lambda'$ and in this case, we define $\xi_{B,A}$ as above.

Let $[A,B]$ be and edge of $\gr$. Of all the edges based at $A$, let $B'$ be such that $[A,B']$ follows $[A,B]$ clockwise around $A$. We define $I_{A\to B}:=[\xi_{A,B},\xi_{A,B'})$. For an edge $[B,A]$, let $A'$ be such that $[B,A']$ follows $[B,A]$ anticlockwise around $B$. In that case, we define $I_{B\to A}:= [\xi_{B,A'},\xi_{B,A})$. The geodesic $\lambda$ and $\lambda'$ will be referred to as the \emph{spectacles defining geodesics} of $I_{A\to B}$ and of $I_{B\to A}$ respectively. We denote this pair of spectacles by $S^f$.

In Section \ref{sec:spectacles}, it is explained how a pair of spectacles $S$ yields $S$-admissible paths. We refer the reader to Lemma 3.5 of \cite{DP14} for the proof that if such paths exist, they must be unique.

As mentionned in Section \ref{sec:spectacles}, it is rather hard to use the geometric interpretation of spectacles to obtain concrete encodings. This is why the encoding in words is important. We recall that by Proposition \ref{prop: admissible_iff_some_inequalities}, to understand $S^f$-admissible paths, we only need to understand what are the limiting words $u_L,u_R,v_L$ and $v_R$ (see Definition \ref{def: extreme-most_words}). However, to compute those words explicitely, we will need some technical results. This will be done in Section \ref{sec: limiting-words}

\section{Monotony of the Busemann functions}\label{sec:monotony_of_bus_fun}

The goal of this section is to show that a semi infinite path $\gamma$ defined by the spectacles $S^f$ does indeed converge to its limit point. As stated in Section \ref{sec:spectacles}, the proof in \cite{DP14} relies on a picture, so we add a this section for completion. Recall that using their Lemma 3.14, this suffices to show that $S^f$ is an sub-accurate pair of spectacles. We now set up some notation and prove some technical lemmas to establish this result.

A semi infinite $S^f$-admissible path $\gamma$ is determined by the data of a vertex 
$x\in V(\gr)$ and a point $\xi\in\del \DD$. In the following, for simplicity, we will suppose that $x$ is a point of type $A$. We will now set up some notation.

\begin{lem}\label{lem:vertex_geod_following_spect_defining_geo_do_not_cross_AB}
	Let $P$ be a face polygon, $e = [A,B]$ one of its edges and $\lambda$ be the spectacles defining geodesic of $I_{A\to B}$. Define another semi infinite geodesic $\alpha$ as follows. Let $A'$ be the vertex of type $A$ coming clockwise after $A$ around $B$; $\alpha$ is the semi infinite geodesic starting at $B$ that goes through $A'$. Then $\alpha$ and $\lambda$ do not cross.
\end{lem}

\begin{figure}
	\centering
	\includegraphics{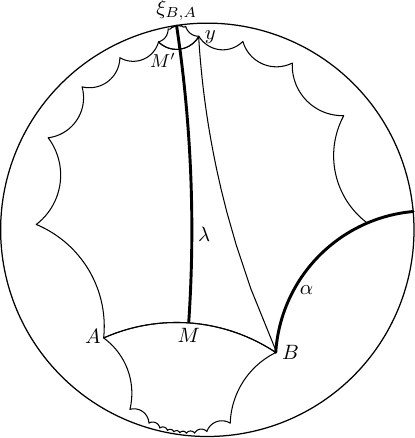} 
	\caption{The geodesics $\alpha$ and $\lambda$ do not cross.
		\label{fig:alpha_lambda_dont_cross}}
\end{figure}

\begin{proof}
	Recall the construction of convex paths done in Lemma \ref{lem: infinite_edge_seq}. Let $e'$ be the opposite vertex to $e$ in $P$, and write $e' = [x,y]$. Denote by $M$ and $M'$ the intersection of $\lambda$ with $e$ and $e'$ respectively. In the course of Lemma \ref{lem: infinite_edge_seq}, we've shown that the quadrilateral formed by $M,M',y,B$ has all its internal angles smaller or equal than $\frac \pi 2$. Hence, $\lambda$ is disjoint from the complete geodesic going through $B$ and $y$. Since the latter crosses $\alpha$ at $B$, $\alpha$ cannot cross $\lambda$ (see Figure \ref{fig:alpha_lambda_dont_cross}).
\end{proof}

In the same order of idea, we have the following Lemma for $I_{B\to A}$-type intervals.

\begin{lem}\label{lem:vertex_geod_following_spect_defining_geo_do_not_cross_BA}
	With the notation above, let $B'$ be the vertex of type $B$ coming anti-clockwise before $B$ around $A$ and let $\alpha'$ be the semi infinite geodesic going through $A$ and $B'$. Then $\alpha'$ does not cross $\lambda'$.
\end{lem}
\begin{proof}
	The proof goes as in Lemma \ref{lem:vertex_geod_following_spect_defining_geo_do_not_cross_AB}
\end{proof}

\begin{figure}[h]
	\centering
	\begin{subfigure}[b]{0.45\textwidth}
		\includegraphics[width=\textwidth]{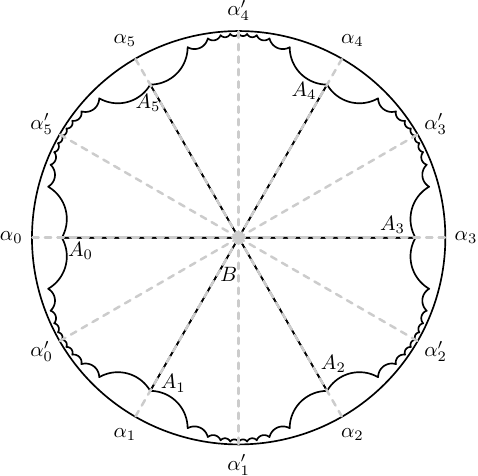}
		\caption{The dashed gray lines represent the $\alpha_j$ and $\alpha_j'$ geodesic segments.}
		\label{fig:spec_def_geod1}
	\end{subfigure}
	\hfill
	\begin{subfigure}[b]{0.45\textwidth}
		\includegraphics[width=\textwidth]{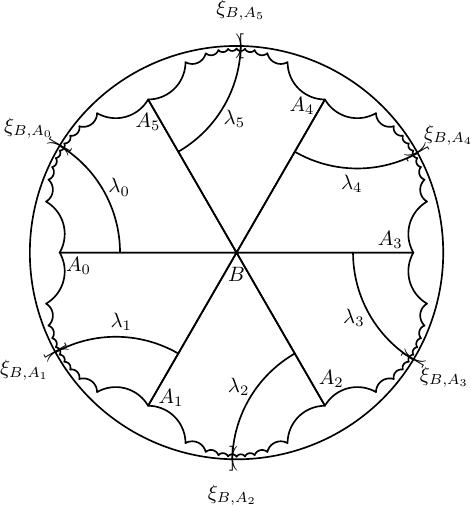}
		\caption{The curved lines indexed by $\lambda_j$ are spectacles defining geodesics.}
		\label{fig:spec_def_geod2}
	\end{subfigure}
	\caption{}
	\label{fig:spec_def_geod}
\end{figure}

Let $(p,q,r)$ be a hyperbolic triplet. Suppose that $\gamma$ is a $S^f$-admissible path that goes through some $A$-type vertex $x$ and a $B$-type vertex $y$ consecutively. Up to conjugating by an isometry, we may assume that $y=0$ and $[x,0]$ aligns with $\RR_{\leq 0}$. Then, there are $q$ choices of points of type $A$ through which $\gamma$ may continue. We label those points anticlockwise $A_0,A_1,...A_{q-1}$, where $A_0 = x$. We define the $2q$ geodesic half segments $\alpha_j:= \left[0, \exp\left(i\pi + j\frac{2i\pi}{q}\right)\right]$ and $\alpha'_j = \left[0,\exp\left(i\left(\pi+\frac{\pi}{2q}\right) + j\frac{2i\pi}{q}\right)\right]$ (see Figure \ref{fig:spec_def_geod1}). Label polygons around $0$ by $P_0,P_1,...,P_{q-1}$ counterclockwise starting with the one on the right of $[A_0,0]$. We write $\lambda_j$ for the spectacles defining geodesics that crosses the segment $\overline{P_j}\cap \overline{P_{j+1}}$, so that $\lambda_j$ is the spectacle defining geodesic for $I_{B\to A_{j+1}}$ (see Figure \ref{fig:spec_def_geod2}).

\begin{prop}\label{prop: Busemann_funct_decreas}
Let $(p,q,r)$ be a hyperbolic triplet. There exists $\varepsilon>0$ such that for all $k\in\ZZ$, the Busemann functions based at $\xi$ satisfy $$B_\xi(x_{k+1}) - B_\xi(x_k) >\varepsilon$$ for all $S^f$-admissible paths $\gamma$, where $(x_k)_{k\in\ZZ}$ is the sequence of points of type $A$ and $\xi$ is the point at $+\infty$ of $\gamma$.
\end{prop}

\begin{proof}
	Let $\gamma$ be a $S^f$-admissible path as above. Fix some base $A$-type vertex of $\gamma$, $x_k$, and move things around so $\gamma$ goes through $x_k = A_0$ followed by $0$ as in the configuration described above.
	
	Since $\lambda_{j+1}$ goes through opposite edges in $P_j$ and $\alpha_j'$ goes through alternating vertices of these edges, they must cross in $P_j$. By Lemma \ref{lem:vertex_geod_following_spect_defining_geo_do_not_cross_AB}, $\lambda_j$ does not cross $\alpha_{j-1}$ for all $j$. By abuse of notation, we will write $\alpha_j$ (resp. $\alpha_j'$) for the point at infinity of the segment $\alpha_j$ (resp. $\alpha_j'$) and $\lambda_j$ for the point at infinity of $\lambda_j$. Since by definition, the $\alpha_j$ and $\alpha_j'$ come in alternating order, we get the cyclic inequality
	
	$$ \alpha_0 < \lambda_0 < \alpha_0' < \alpha_1 < \lambda_1 < \alpha'_1 < ... < \alpha_{r-1} < \lambda_{r-1} < \alpha_{r-1}' < \alpha_0.$$
	
	We also need to consider the spectacles defining geodesic $\beta$ imposed by the edge $[A_0,0]$. Remark that, $\beta$ must cross $\alpha_0'$ and $\lambda_0$ in $P_0$, but by Lemma \ref{lem:vertex_geod_following_spect_defining_geo_do_not_cross_BA}, it does not cross $\alpha_1$. As before, if we write $\beta$ for the endpoint of $\beta$, this gives us the cyclic inequality $$\alpha_0<\lambda_0<\alpha_0'<\beta<\alpha_1 <\alpha_0.$$
	
	Combining the two previous inequality yields 
	\begin{equation}\label{ineq:in_circle}
		\alpha_0 < \lambda_0 < \alpha_0' < \beta< \alpha_1 < \lambda_1 < \alpha'_1 < ... < \alpha_{r-1} < \lambda_{r-1} < \alpha_{r-1}' < \alpha_0,
	\end{equation}
	which is illustrated in Figure \ref{fig:all_curves}.
	\begin{figure}[h]
		\centering
		\includegraphics{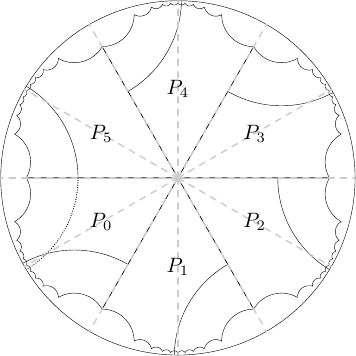}
		\caption{The cyclic inequality \eqref{ineq:in_circle} says that endpoints of the $\lambda_j$, $\alpha_j$ and $\alpha_j'$ come in this order. The dotted line represents $\beta$.}
		\label{fig:all_curves}
	\end{figure}
	
	Denote by $\xi$ the positive endpoint of $\gamma$. Recall that $\gamma$ goes through $A_0,0,A_j$ if and only if $\lambda_{j}<\xi \leq\lambda_{j+1}$. Also, since $\gamma$ first goes through $A_0$ and $0$, we have $\beta < \xi \leq \lambda_{q-1}$. So for $\gamma$ to go through $A_0$ again, we would need to have $\lambda_{q-1} < \xi \leq \lambda_0$, but this is impossible since by \eqref{ineq:in_circle}, $\lambda_{q-1} < \lambda_0 < \beta < \lambda_{q-1}$. Furthermore, $\gamma$ goes through $A_0,0,A_1$ if and only if $ \beta <\xi \leq \lambda_1$.
	
	Consider $$\mathcal{B}_j:=\{\zeta\in S^1:B_\zeta(A_{k+1}) - B_\zeta(A_k)>0\}.$$
	
	To show that Busemann functions are increasing, we need to show that, for $2\leq j\leq q-1$, if  $\lambda_{j-1} < \xi \leq \lambda_{j}$, then $\xi\in\mathcal{B}_j$ and if $\beta < \xi \leq \lambda_1$, then $\xi\in\mathcal{B}_1$.
	
	If $r$ is even, then $$\xi\in\mathcal{B}_j\iff
	\begin{cases}
		\alpha'_n < \xi < \alpha'_{n+\frac r 2},\quad \text{if $j = 2n+1$ for some $0\leq n\leq \frac r 2 -1$};\\
		\alpha_n < \xi <\alpha_{n+\frac r 2},\quad \text{if $j = 2n$ for some $1\leq n\leq \frac r 2 -1$}.
	\end{cases}$$
	
	First suppose that $j = 2n+1$ and $1\leq n \leq \frac r 2 -1$. Then in particular, $n\leq 2n-1$ and $2n+1\leq n+\frac r 2$, so that using inequality \eqref{ineq:in_circle}, we have
	$$\alpha'_n\leq \alpha'_{2n-1} < \alpha_{2n} < \lambda_{2n} < \lambda_{2n+1} < \alpha_{2n+1}'\leq \alpha_{n +\frac r 2}'.$$
	
	Now, suppose that $j = 2n$ and $1\leq n \leq \frac r 2 -1$. Then in particular, $n\leq 2n-1$ and $2n+1 \leq n +\frac r 2$. Hence, $$\alpha_n\leq \alpha_{2n-1} < \lambda_{2n-1} < \lambda_{2n} < \alpha_{2n+1}\leq \alpha_{n+\frac r 2}.$$
	
	We now have to show that if $\beta < \xi \leq \lambda_1$, then $\alpha_0' < \xi < \alpha_{\frac r 2}'$. Which is clear by Inequality \eqref{ineq:in_circle}. Hence, in all cases, Busemann functions are strictly decreasing. 
	
	If $r$ is odd, then $$ \xi\in \mathcal{B}_j \iff 
	\begin{cases}
		\alpha_n' < \xi < \alpha_{n+\frac {r+1} {2}},\quad \text{if $j = 2n+1$ for some $0\leq n \leq \frac {r-3}{2}$;}\\
		\alpha_n < \xi < \alpha_{n+\frac {r-1}{2}},\quad\text{if $j=2n$ for some $1\leq n \leq \frac{r-1}{2}$.}
	\end{cases}$$
	
	The proof is the same as above.
	
	Hence, for any $S^f$-admissible path $\gamma$, and for every $k\in \ZZ$, $B_\xi(x_{k+1}) - B_\xi(x_k) >0$. Since this difference attains its extremas at boundary points of the intervals $I_{0\to A_j}$ and there are only finitely many such endpoints, we can take $\varepsilon := \min B_\xi(x_{k+1}) - B_\xi(x_k)$, where the minimum is taken over endpoints of the interval corresponding to the choice $x_{k+1}$ and over all $q-1$ such choices.
\end{proof}

\section{Limiting words of $S^f$}\label{sec: limiting-words}

We are now ready to determine the codes of the limiting words for the pair of spectacles $S^f$. 

\begin{prop}\label{prop: extreme-most words for the spectacles Sf}
	Suppose that $p\geq 3$. The limiting words are given by the following table
	
	\begin{center}\label{table: codes}
		\begin{tabular}{c|c|c}
			& $u_L$& $v_R$\\
			\hline\hline
			$p\geq 3$ & &\\
			\hline
			$r$ odd & $((a^{p-1}b)^{\frac{r-3}{2}}a^{p-1}b^2)^\infty$ & $(b^{q-1}a)^{\frac{r-3}{2}} b^{q-1}a^2)^\infty$\\
			\hline
			$r$ even  &$((a^{p-1}b)^{\frac{r-2}{2}}a^{p-2}(ba^{p-1})^{\frac{r-2}{2}}b^2 )^\infty$ & $((b^{q-1}a)^{\frac{r-2}{2}}b^{q-2}(ab^{q-1})^{\frac{r-2}{2}}a^2)^\infty$\\
			\hline
		\end{tabular}
	\end{center}
	
\end{prop}

\begin{proof}
	We will do the proof for the case where $r$ is odd and for the computation of $u_L = w(\gamma^{S^f}_{B_0\to \xi_{B_0,A_0}})$.
	
	Let $P$ be a face polygon and label its edges from $0$ to $2r-1$, where each edge has the orientation induced by the orientation of $P$. Label vertices of $P$ from $0^+$ to $(2r-1)^+$ in a way that $(2r-1)^+ = A_0$ and $0^+ = B_0$. 
	
	Our first claim is that the first $r-1$ edges of $\gamma_{B_0\to\xi_{B_0,B_0}}^{S^f}$ are $1,2,3,...,r-1$. To prove this, we will show that $\xi_{B_0,A_0}\in I_{0^+\to 1^+}\cap I_{1^+\to 2^+}\cap ... \cap I_{(r-2)^+\to (r-1)^+}$. Throughout the proof, we will write $\lambda$ for the spectacle defining geodesic of $I_{B_0\to A_0}$.
	
	Let $1\leq k\leq r-1$. First suppose that $k^+$ is a point of type $A$. Let $B$ be the point of type $B$ coming clockwise after $(k+1)^+$ around $k^+$. We will write $\alpha$ for the spectacle defining geodesic of $I_{k^+\to (k+1)^+}$ and $\alpha'$ as the spectacle defining geodesic of $I_{k^+\to B}$. As in the proof of Proposition \ref{prop: Busemann_funct_decreas}, we will identify half geodesics with their endpoints onf $\del \DD$. To show the claim, we will prove that $\alpha' < \lambda < \alpha$. First, remark that by Proposition \ref{prop: Busemann_funct_decreas}, the half-geodesics $\alpha$ and $\lambda$ cross in $P$. Recall that by definition of spectacles defining geodesics, and by Lemma \ref{lem: infinite_edge_seq}, the half geodesic $\alpha$ crosses sides $k+1$ and $k+1+r$. Let $\beta$ be the half geodesic starting at $(k+1+r)^+$ and going through antipodal vertices in $P$. As in the proof of Proposition \ref{prop: Busemann_funct_decreas}, we have that $\beta < \lambda < \alpha$. Furthermore, if we consider $\beta$ as a half geodesic starting in $(k+1)^+$, we obtain the inequality $\alpha' < \beta < \lambda$. Those two inequalities imply that $\alpha' < \lambda < \alpha$. The proof is similar when $k^+$ is of type $B$. Since the edge $[0^+,0^-]$ my be pushed to $[r^-, r^+]$ via a transformation fo $G_{p,q,r}$, the picture near $B_0$ is the same as the one near $r^-$. Hence the path is periodic, which gives the code $u_L = ((a^{p-1}b)^{\frac{r-3}{2}}ab^2)^\infty$. The proof is similar for the other words.
	\end{proof}
	
	As stated in Section \ref{sec:spectacles}, this property is very important. To give a concrete enumeration of closed geodesics, we now only have to enumerate some words in an alphabet that satisfy come inequalities. We also get as an easy corollary an important property of $S^f$-admissible paths.

	\begin{cor}\label{cor: paths_take_less_than_r_sides}
		Let $\gamma$ be a $S^f$-admissible path. Then for each face polygon, $\gamma$ goes through less than half of its sides.
	\end{cor}
	
	\begin{proof}
		If $\gamma$ took more than half of the sides of a face polygon, its code wouldn't satisfy the inequalities prescribed by Proposition \ref{prop: admissible_iff_some_inequalities}.
	\end{proof}

\section{Construction of paths framing words}\label{sec:construction_of_paths_framing_words}

The work that has been done so far gives us some language to talk about $S^f$-admissible paths, but not much has been said about the interplay between a geodesic and its associated $S^f$-admissible path. In this section, we will show that they are geometrically close one to another.

\begin{lem}\label{lem:ad_path_are_straight}
	Let $\gamma$ be an $S^f$-admissible path joining $\eta$ to $\xi$, where $\eta,\xi\in\del\DD$. Then $\gamma$ is an embedded proper path.
\end{lem}

\begin{proof}
	Since all of $\gamma$'s edges are embedded, we only need to show that $\gamma$ has no cycle. For this, write $\gamma = (...,A_{-2},B_{-1},A_0,B_1,A_2,...)$. By Proposition \ref{prop: Busemann_funct_decreas}, $A$-types vertices of $\gamma$ are all different. Hence, if a cycle existed in $\gamma$, it would emerge from indices $k_1,k_2\in\ZZ$ for which $B_{2k_1+1} = B_{2k_2+1}$. However, the choice of point of type $A$ that follows $B_{2k_i+1}$ is determined by the endpoint $\gamma_{+\infty}$. This implies that $A_{2k_1+2} = A_{2k_2+2}$, which is impossible.
\end{proof}

\begin{figure}
	\centering
	\includegraphics{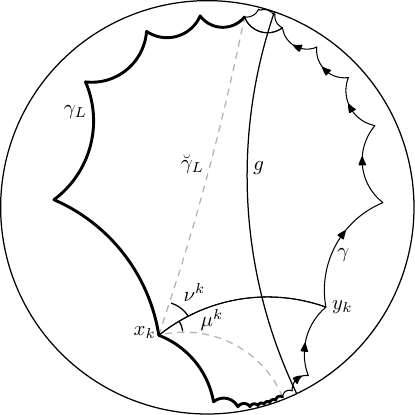} 
	\caption{The path in bold represents $\gamma_L$, the dashed one, $\breve{\gamma}_L$ and the one with arrows on it, $\gamma$. The line indexed by $g$ is the geodesic that defines $\gamma$. 
		\label{fig:the_three_paths}}
\end{figure}

\begin{prop}\label{prop:existence_of_cvx_paths}
	Let $\gamma$ be a $S^f$-admissible path. There exists embedded proper paths $\breve{\gamma}_L$ and $\breve{\gamma}_R$ (which do not lie on $\gr$) that are respectively right and left convex. These paths satisfy $$ \breve{\gamma}_L\prec\gamma\prec\breve{\gamma}_R.$$
	
	In particular, there exists a unique geodesic $g'$ satisfying $\breve{\gamma}_L\leq g'\leq \breve{\gamma}_R$. If $\gamma$ if associated to a geodesic $g$, necessarily, $g =g'$.
\end{prop}
\begin{proof}
	We will first construct $\breve{\gamma}_L$; its right-convexity will be a consequence of Lemma \ref{lem:generic_case_smaller_than_pi2}. The same will be true for $\breve{\gamma}_R$. The desired inequalities will be self-evident from the construction.
	
	The edges $[x,y]$ in the region $\overline{L_\gamma}$ that are not edges of $\gamma$ either have one or zero vertex in common with $\gamma$. Let $\mathcal{EL}_\gamma$ be the set of such edges that have one vertex in common with $\gamma$ and write $\gamma = (\gamma_j)_{j\in\ZZ}$.
	
	For two edges $e_1,e_2\in\mathcal{EL}_\gamma$, there exists $j_1$ and $j_2$ such that $\gamma_{j_i}$is the vertex that $\gamma$ and $e_i$ have in common. We will say that $e_1\leq e_2$ if $j_1\leq j_2$. Hence, we can write $$\mathcal{EL}_\gamma = (e_k)_{k\in\ZZ}$$ as a totally ordered set. 
	
	For $k\in \ZZ$, we write $y_k$ for the vertex that $e_k$ shares with $\gamma$ and $x_k$ for the other vertex. We define $$\breve{\gamma}_L:=\bigcup_{k\in\ZZ} [x_{k},x_{k+1}].$$
	
	To show that $\breve{\gamma}_L$ is right convex, it suffices to show that for every $k\in\ZZ$, the angle $\theta^k:=\angle x_{k-1},x_k,x_{k+1}$ is acute. If we let $\mu^k:=\angle x_{k-1},x_k,y_k$ and $\nu^k:=\angle y_k, x_k,x_{k+1}$, then $\theta^k = \nu^k+\mu^k$. Recall that by Corollary \ref{cor: paths_take_less_than_r_sides}, a $S^f$-admissible path may go to at most $r$ sides of any face polygon. Hence, for every $k\in\ZZ$, $\nu^k$ or $\mu^k$ is at most the angle corresponding to this maximal configuration. For hyperbolic triplets $(p,q,r)$ with $p\geq 3$ and $(p,q,r)\not = (3,3,4)$ or $(3,4,4)$, we've shown in Lemma \ref{lem:generic_case_smaller_than_pi2} that both angles $\mu^k$ and $\nu^k$ are acute, so that $\theta^k$ is convex (see Figure \ref{fig:the_three_paths}).
	
	In the cases $(3,3,4)$ and $(3,4,4)$, we must do explicit computations. Let $y_k$ be some vertex of $\gamma$ and denote by $z_k$ its following vertex in $\gamma$. We will write $P$ for the polygon on the left of $[y_k,z_k]$ and $P'$ for the polygon on its right. Since $r = 4$ is even, by Corollary \ref{cor: paths_take_less_than_r_sides} the maximal configuration in this case corresponds to $\gamma$ going through $4$ sides of $P'$ and $3$ sides of $P$. Direct computations show that in both cases, $\theta^k < \pi$. 
	
	To construct the paths $\breve{\gamma}_R$, we do the same procedure considering the set $\mathcal{ER}_\gamma$ of edges that have only one vertex in common with $\gamma$ laying in $\overline{R_\gamma}$.
	
	Since $\gamma$ is embedded, $\breve{\gamma}_L$ and $\breve{\gamma}_R$ are also embedded. By construction, we have $$\breve{\gamma}_L \prec \gamma \prec \breve{\gamma}_R.$$
\end{proof}

We will also need to see how these convex paths interact with graph paths near a $S^f$-admissible path $\gamma$.

\begin{defn}\label{def:ajd_interm_polygons}
	\begin{itemize}
		\item A face polygon $P$ of $\gr$ is a \emph{neighbouring polygon of $\gamma$} if $\overline{P}\cap \gamma\not = \varnothing$.\\
		\item For a neighbouring polygon $P$, either $P\cap\gamma$ contains an edge or it does not. We will refer to the former as an \emph{adjacent polygon to $\gamma$} and to the latter as an \emph{intermediate polygon to $\gamma$}.
	\end{itemize}

\end{defn}
By Lemma \ref{lem:ad_path_are_straight} and the discussion in Section \ref{sec:the_graph}, we are also able to talk about the left and right version of polygons defined above in reference to $S^f$-admissible paths. 

\begin{defn}
	Let $\gamma$ be a $S^f$-admissible path. The \emph{band on the left of $\gamma$} is $$\joliL_\gamma:=\{\text{neighbouring polygons of $\gamma$ on its left}\}.$$ 
	
	Analogously, we define $\joliR_{\gamma}$, \emph{the band on the right of $\gamma$}.
\end{defn}

\begin{rem}
	As though we've defined $\joliL_\gamma$ and $\joliR_\gamma$ as sets, most of the time, we will want to consider them as totally ordered sets. Indeed, the orientation on $\gamma$ lets us think of a, say, left-neighbouring polygon $P$ coming before another left-neighbouring polygon $Q$.
\end{rem}

If $\gamma$ is a $S^f$-admissible path, recall that by Corollary \ref{cor: paths_take_less_than_r_sides}, $\gamma$ goes through at most $r$ sides of each of its left-neighbouring (resp. right-neighbouring) polygons. Hence, for such a polygon $P$, some of its edges have disjoint closure with $\gamma$. If we join these edges together for each left-neighbouring (resp. right-neighbouring) polygon, we obtain a path $\gamma_L$ (resp. $\gamma_R$).

\begin{prop}\label{prop:ineq_on_5_paths}
	The paths $\gamma_L$ and $\gamma_R$ are embedded proper paths. Furthermore, we have that $$\gamma_L\ll\breve{\gamma}_L \prec \gamma \prec \breve{\gamma}_R \ll\gamma_R.$$
\end{prop}
\begin{proof}
	We can write $$\gamma_L = \bigcup_{\text{left-neighbouring polygons $P$}}\gamma_L\cap P.$$ Each segment $\gamma_L\cap P$ is embedded as the boundary of a polygon. Since each left-neighbouring polygon $P$ is distinct, and they converge to the point at $+\infty$ of $\gamma$, the segments $\gamma\cap P$ glue together in a global embedding. The same is true for $\gamma_R$. 
	
	The inequality $$\gamma_L\ll\breve{\gamma}_L \prec \gamma \prec \breve{\gamma}_R \ll\gamma_R$$
	is direct from the construction of $\gamma_L$ and $\gamma_R$.
\end{proof}

The next lemma shows that the paths $\gamma_L$ and $\gamma_R$ behave sort of like the floor and ceiling functions among some subset of $\Sigma_{\eta,\xi}$.

\begin{lem}\label{lem:partie_entière}
	Let $\eta \not= \xi \in\del\DD$, $\alpha\in\mathcal{P}(\gr)\cap\Sigma_{\eta,\xi}$ and $\gamma$ be a $S^f$-admissible path joining $\eta$ to $\xi$. If $\gamma\prec\alpha$, then $\gamma_R\leq\alpha$. In a similar fashion, if $\alpha\prec\gamma$, then $\alpha\leq \gamma_L$.
\end{lem}

\begin{proof}
	We will only show the first inequality; the other one will fall under the same logic. Suppose that $\alpha \not \geq \gamma_R$. Then there exists a point $x$ of $\alpha$ such that $x\in L_{\gamma_R}$. We will first show that we can suppose that this point can be taken as a vertex of $\alpha$. Indeed, if the point $x$ is in the interior of an edge $e = [y,z]$, since $x\in L_{\gamma_R}$, it follows that the open edge $(y,z)\subseteq L_{\gamma_R}$. Hence, if we had $y,z\in \overline{R_{\gamma_R}}$, since they are limit points of $(y,z)$, we would need to have that $y,z\in\del\overline{R_{\gamma_R}} = \gamma_R$. But this forces $[y,z]$ to be an edge of $\gamma_R$ which is impossible since this would imply that $x\in\gamma_r\subseteq\overline{R_{\gamma_R}}$. So, without loss of generality, we will assume that $y\in L_{\gamma_R}$. Since $\gamma\prec \alpha$, we conclude that $y\in R_\gamma\cap L_{\gamma_R}$. This is impossible because there are no vertex of $\gr$ in $R_\gamma \cap L_{\gamma_R} = \varnothing$ by construction.
\end{proof}

\section{Almost-unicity of the representation in words}\label{sec:almost_unicity}

In \cite{DP14}, it is shown how the uniqueness of representation in words fails. Here, we will give a different proof of this fact.

\begin{thm}\label{thm:unicity}
	Let $\gamma$ and $\omega$ be periodic $S^f$-admissible paths such that ${\gamma}_{\infty} = { \omega}_\infty$ and ${\gamma}_{-\infty} = { \omega}_{-\infty}$. Then if $\gamma \not = \omega$, we have $\{w(\gamma),w(\omega)\} = \{w_L,w_R\}$, where $w_L$ and $w_R$ are the words defined in Section \ref{sec:intro}.
\end{thm}

The latter theorem is a consequence of the results of the previous sections and the following lemmas.

\begin{lem}\label{lem:containmentimplieseq}
	Suppose that $\gamma$ is a $S^f$-admissible path and that $\alpha\in\mathcal{P}(\gr)\cap\Sigma_{\gamma_{-\infty},\gamma_\infty}$. If $\breve{\gamma}_L \ll \alpha \leq \gamma$ or $\gamma\leq \alpha \ll \breve{\gamma}_R$, then $\alpha = \gamma$.
\end{lem}

\begin{proof}
	The proof goes by contradiction. We will suppose that $\breve{\gamma}_L \ll \alpha\leq \gamma$; if the other condition is satisfied, the proof is a consequence of an analog reasoning. Write $\alpha = (\alpha_j)_{j\in\ZZ}$ and $\gamma = (\gamma_j)_{j\in \ZZ}$. Suppose that $k$ is such that $\alpha$ never goes through $\gamma_k$. Denote by $\breve{\mathcal{S}}_L$ the closed strip contained in between $\breve{\gamma}_L$ and $\gamma$ and let $G$ be the sub-graph induced by vertices in $\breve{\mathcal{S}}_L$. If we remove $\gamma_k$ and its neighbouring edges, then $G = G_{-\infty}\sqcup G_\infty$, where $\overline{G_{\pm\infty}}\ni \gamma_{\pm\infty}$. Since $\alpha$ is connected, only one of $\alpha \subseteq G_\infty$ or $\alpha\subseteq G_{-\infty}$ can hold; in either case, $\alpha$ will not be able to join its other endpoint. 
\end{proof}

\begin{proof}(of Theorem \ref{thm:unicity})
	First remark that if $\gamma$ and $\omega$ had a point $P$ in common, they would agree after $P$ by the unicity of semi-infinite admissible paths. Hence, in the following, we will suppose that $\gamma$ and $\omega$ are disjoint. Without loss of generality, we assume that $\gamma \prec \omega$. By Proposition \ref{prop:ineq_on_5_paths}, there exists convex paths satisfying \begin{equation}\label{eq:ineq_gamma}
		\gamma_L\ll \breve{\gamma}_L\prec\gamma\prec\breve{\gamma}_R\ll\gamma_R
	\end{equation}
	and
	\begin{equation}\label{ineq_omega}
		\omega_L\ll\breve{\omega}_L\prec\omega\prec\breve{\omega}_R\ll\omega.
	\end{equation}
	
	Because $\gamma$ and $\omega$ are embedded proper graph paths satisfying $\gamma\prec\omega$, by Lemma \ref{lem:partie_entière}, we have 
	\begin{equation}\label{ineq_int1}
		\gamma_R\leq\omega
	\end{equation}
	and
	\begin{equation}\label{ineq_int2}
		\gamma\leq \omega_L.
	\end{equation}
	
	Furthermore, if $g_\gamma$ and $g_\omega$ are the geodesics given by Proposition \ref{prop:existence_of_cvx_paths}, since $\gamma$ and $\omega$ converge to the same point at $\pm\infty$, we have $g:=g_\gamma = g_\omega$ and
	\begin{equation}\label{geod_gamma}
		\breve{\gamma}_L\leq g\leq \breve{\gamma}_R
	\end{equation} 
	and
	\begin{equation}\label{geod_omega}
		\breve{\omega}_L\leq g \leq \breve{\omega}_R.
	\end{equation}
	
	Combining \eqref{geod_gamma} and \eqref{geod_omega} together yields
	
	\begin{equation}\label{comp2}
		\breve{\omega}_L\leq \breve{\gamma}_R
	\end{equation}
	                                                                           
	Now, if we use \eqref{ineq_int2} with \eqref{ineq_omega} and \eqref{comp2}, we obtain
	\begin{equation*}\label{eq:bounds1}                        
		\gamma\leq\omega_L\ll\breve{\omega}_L\leq\breve{\gamma}_R.
	\end{equation*}
	By Lemma \ref{lem:containmentimplieseq}, it follows that $\omega_L = \gamma$.
	
	In a similar way, using \eqref{comp2}, \eqref{eq:ineq_gamma} and \eqref{ineq_int1}, we obtain
	\begin{equation*}
		\breve{\omega}_L\leq\breve{\gamma}_R\ll\gamma_R\leq\omega,
	\end{equation*}
	which, again, implies that $\gamma_R = \omega$.
	
	Hence, for such a situation to happen, the two paths $\gamma$ and $\omega$ must be neighbours of each other.
	
	The argument at the end of Lemma 3.15 in \cite{DP14} shows that this is only possible when $\gamma$ and $\omega$ have the codes $w_L$ and $w_R$.
\end{proof}

\section{A concrete algorithm to compute the Length Spectrum of triangular surfaces}\label{sec:main_thm}

Let $(p,q,r)$ be a hyperbolic triplet with $p\geq 3$. The results of Section \ref{sec:almost_unicity} give a unique encoding for every (but one) geodesic of $X_{p,q,r}$. The first thing to do is to understand what this encoding means geometrically. Recall that if $w$ is the code of a $S^f$-admissible path $\gamma$, then letters $a^e$ correspond to some amount of turning around an $A$-type vertex. If $e = 1$ or $p-1$, then the path $\gamma$ continues to follow some face polygon. Else, $2\leq e\leq p-2$ and $\gamma$ shifts from one face polygon to another. The same goes for letters $b^f$ whether $b=1$, $q-1$ or $2\leq b \leq q-2$. We will refer to  the first type of letters as \emph{switches}.

\begin{figure}[h]
	\centering
	\begin{subfigure}[b]{0.45\textwidth}
		\includegraphics[width=\textwidth]{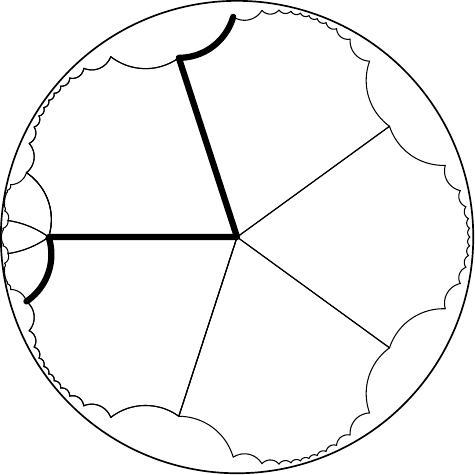}
		\caption{The word $(ab)^\infty$ is made of one zigzag.}
		\label{fig:zigzag1}
	\end{subfigure}
	\hfill
	\begin{subfigure}[b]{0.45\textwidth}
		\includegraphics[width=\textwidth]{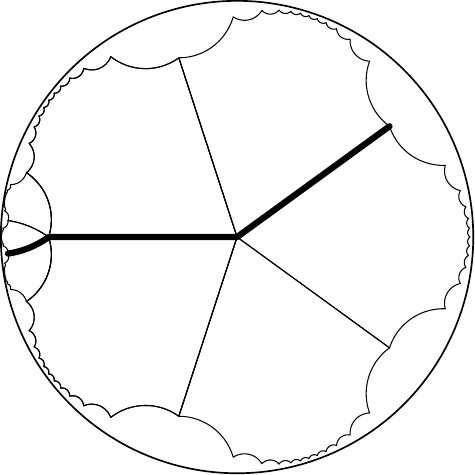}
		\caption{The word $(a^2b^2)^\infty$ decomposes as $...(a^2)()(b^2)()(a^2)()...$, where $()$ denotes the empty zigzag. Geometrically, the empty zigzag corresponds to $\gamma$ going through only on side of an adjacent polygon.}
		\label{fig:zigzag2}
	\end{subfigure}
	\caption{}
	\label{fig:zigzag}
\end{figure}

Accordingly to the discussion above, every word factors into switches and what we call \emph{zigzags}, i.e. subwords that do not contain letters $a^e$ or $b^f$ for $2\leq e\leq p-2$ and $2\leq f \leq q-2$. Remark that a word may contain no switch (e.g. $w=(ab)^\infty$ on Figure \ref{fig:zigzag1}). However, we will want to think of words like $(a^2b^2)^\infty$ as having infinitely many zigzags, meaning that we consider the empty word as a zigzag (see Figure \ref{fig:zigzag2}).

The goal of the following discussion is to show that for a geodesic $g$, the encoding of the $S^f$-admissible path $\gamma(g)$ given in Section \ref{sec:almost_unicity} implies lower bounds for the geometrical length of $g$. To do so, we will use the factorization into switches and zigzags introduced above.

\begin{figure}[h!]
	\centering
	\begin{subfigure}[b]{0.45\textwidth}
		\includegraphics[width=\textwidth]{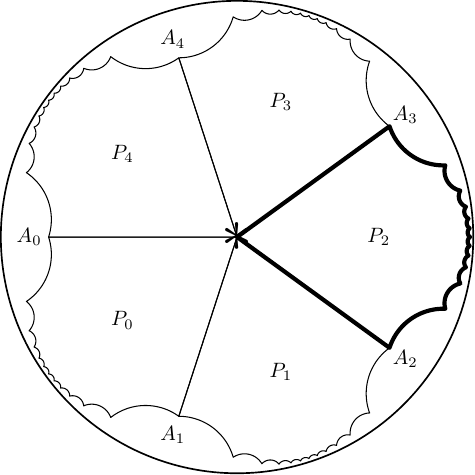}
		\caption{Polygons around the $B$-type vertex $0$.}
		\label{fig:geodesic_set1}
	\end{subfigure}
	\hfill
	\begin{subfigure}[b]{0.45\textwidth}
		\includegraphics[width=\textwidth]{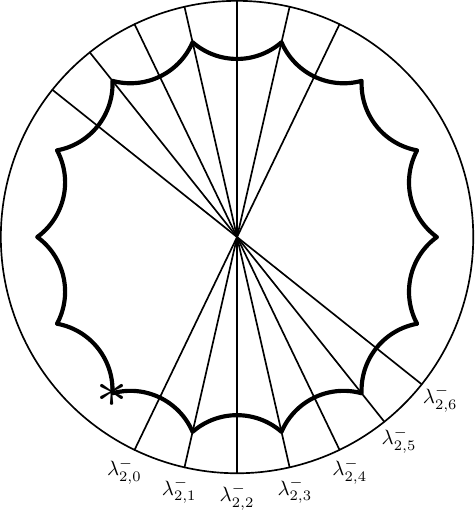}
		\caption{The geodesics $\lambda_{j,k}^-$ for the polygon $P_2$, in bold in Figure \ref{fig:geodesic_set1}.}
		\label{fig:geodesic_set2}
	\end{subfigure}\\
	
	\begin{subfigure}[b]{0.45\textwidth}
		\includegraphics[width=\textwidth]{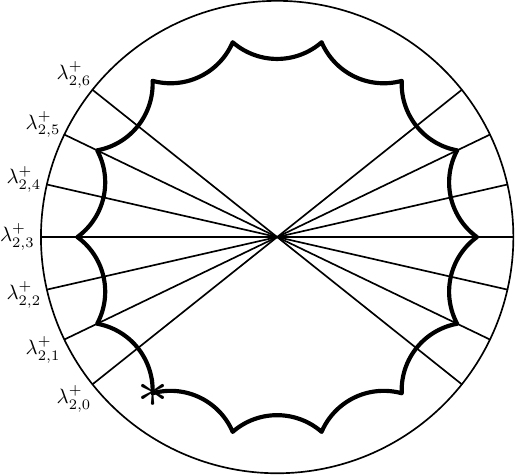}
		\caption{The geodesics $\lambda_{j,k}^+$ for the polygon $P_2$, in bold in Figure \ref{fig:geodesic_set1}.}
		\label{fig:geodesic_set3}
	\end{subfigure}
	\hfill
	\begin{subfigure}[b]{0.45\textwidth}
		\includegraphics[width=\textwidth]{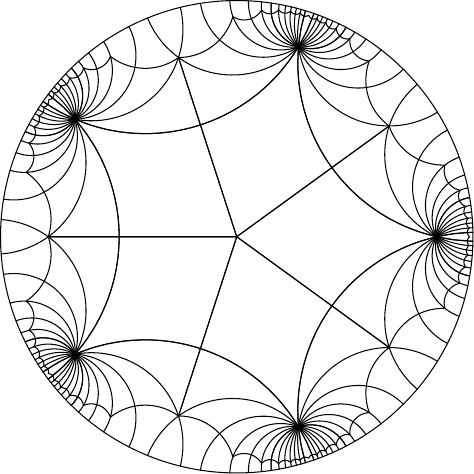}
		\caption{The set of geodesics $\lambda_{j,k}^\pm$ around $0$.}
		\label{fig:geodesic_set4}
	\end{subfigure}
	\caption{}
	\label{fig:geodesic_set}
\end{figure}

For now, let $\gamma = \gamma(g)$ be the $S^f$-admissible path of some geodesic $g$ and $w = w(\gamma)$ be its encoding. Let $\gamma_{m-1},\gamma_{m},\gamma_{m+1}$ be successively $A$-type, $B$-type and $A$-type vertices. Up to conjugation by an isometry, we will suppose that $\gamma_m = 0$ and $\gamma_{m-1}$ is aligned with $\RR_{<0}$ in $\DD$ (see Figure \ref{fig:geodesic_set1}). Label polygons around $0$ as it was done in the discussion preceding Proposition \ref{prop: Busemann_funct_decreas}. For each polygon $P_j$, label the sides of $P_j$ counterclockwise from $0$ to $2r-1$ starting with $[0,A_j]$ and its vertices from $0^+$ to $2r-1^+$ starting from $A_j$. Then define, for $0\leq k \leq r-1$, $2r$ geodesics as follows. 

\begin{itemize}
	\item If $k$ is even, $\lambda_{j,k}^-$ is the geodesic coming from Lemma \ref{lem: infinite_edge_seq} that goes through edges $\frac k 2$ and $\frac{k}{2} + r$ of $P_j$ and $\lambda_{j,k}^+$ is the geodesic coming from Lemma \ref{lem: infinite_edge_seq} going through edges $2r-1 - \frac k 2$ and $2r-1 - \frac k 2 -\frac r 2$ of $P_j$, see Figures \ref{fig:geodesic_set2} and \ref{fig:geodesic_set3}.\\
	\item If $k$ is odd, then $\lambda_{j,k}^-$ is the geodesic going through vertices $\left(\frac{k-1}{2}\right)^+$ and $\left(\frac{k-1}{2} + \frac r 2\right)^+$ and $\lambda_{j,k}^+$ is the geodesic going through vertices $\left(2r-2 -\frac{k-1}{2} \right)^+$ and $\left(2r-2 -\frac{k-1}{2} - \frac{r}{2}\right)^+$, see Figures \ref{fig:geodesic_set2} and \ref{fig:geodesic_set3}.
\end{itemize}
Some geodesics $\lambda_{j,k}^\pm$ are illustrated in Figure \ref{fig:geodesic_set4} for the other polygons around $0$.

In the same fashion, we consider the case where $\gamma_{m-1},\gamma_m$ and $\gamma_{m+1}$ are successively $B$-type, $A$-type and $B$-type vertices and make similar definitions.

Our next goal is to show that some of the geodesics that have been defined above do not cross. This is done more easily by using the following lemma.

\begin{lem}\label{lem:geod_dont_cross}
	Let $(p,q,r)$ be a hyperbolic triplet and $P$, a face polygon. Suppose that $0^+$ is a vertex of type $A$. Let $\alpha$ be the complete geodesic going through $0^+$ and $1^+$ and $\beta$ be the complete geodesic crossing antipodal vertices in $P$ that meets $(-1)^+$ (and $(-1+r)^+$). Then $\alpha$ and $\beta$ do not cross.
\end{lem}

\begin{figure}[h]
	\centering
	\includegraphics{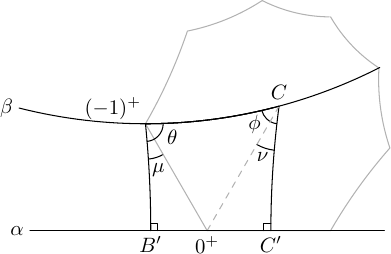}
	\caption{The quadrilateral $\mathcal{Q}$ for $p = 3$.}
	\label{fig:quadrilateral}
\end{figure}

\begin{proof}
	First suppose that $p\geq 4$. In this case, the angles at vertices of $P$ are acute. The quadrilateral formed by $\mathcal{Q} := (-1)^+,0^+,1^+$ and $(-1+r)^+$ has acute angles. Hence, the perpendicular that joins $[0^+, 1^+]$ and $[(-1)^+, (-1+r)^+]$ does so inside $\mathcal{Q}$. Since those two sides are sub-segments of $\alpha$ and $\beta$ respectively, the perpendicular we have exhibited is the common perpendicular to $\alpha$ and $\beta$.
	
	Now, suppose that $p = 3$. Remark that $\beta$ goes through the center $C$ of $P$. Consider the projection $C'$ of $C$ onto $\alpha$ and $B'$ of $(-1)^+$ onto $\alpha$ (see Figure \ref{fig:quadrilateral}). Our goal is to show that the quadrilateral $\mathcal{Q}$ passing through $(-1)^+, B', C'$ and $C$ has acute angles. Let $\mu:=\angle(0^+, (-1)^+, B')$, $\theta:= \angle (B', (-1)^+, C)$ and $\phi := \angle(C', C, (-1)^+)$, we will show that $\theta,\phi \leq \frac \pi 2$.
	
	Since $\Delta(B', (-1)^+, 0^+)$ is a right-angled triangle, we have, using the formulas from the fundamental triangles, that 
	\begin{align*}
		\cot\mu &=\frac{\cos\frac\pi r + \cos \frac \pi 3 \cos \frac \pi q}{\sin \frac \pi 3 \sin \frac \pi q} \tan\frac \pi 3 \\
		&=\frac{\cos \frac \pi r + \cos \frac \pi 3 \cos \frac \pi 3}{\cos \frac \pi 3 \sin \frac \pi q}\\
		&\geq\frac{\cos \frac \pi 4 +\cos \frac \pi 3 \cos \frac \pi 3}{\cos \frac \pi 3 \sin \frac \pi 3}.
	\end{align*}
	Solving for $\mu$ gives that $\theta = \mu + \frac \pi q \lessapprox 1.472... \leq \frac \pi 2$.
	
	For $\phi$, remark that $\phi = \angle(C',C,(-1)^+)\leq \angle(1^+,C,(-1)^+) = \frac  {2\pi} r \leq \frac \pi 2 $.
	
	Hence all internal angles of $\mathcal{Q}$ are acute and the common perpendicular between $\alpha$ and $\beta$ is realized in $\mathcal{Q}$.
\end{proof}

We have an analog lemma for the case where $0^+$ is of type $B$.

\begin{lem}\label{lem:geod_dont_cross2}
		Let $(p,q,r)$ be a hyperbolic triplet and $P$, a face polygon. Suppose that $0^+$ is an edge of type $B$. Let $\alpha$ be the complete geodesic going through $0^+$ and $1^+$ and $\beta$ be the geodesic that crosses antipodal vertices in $P$ that meets $(-1)^+$ (and $(-1+r)^+$). Then $\alpha$ and $\beta$ de not cross.
\end{lem}
\begin{proof}
	The proof is the same as in Lemma \ref{lem:geod_dont_cross} for the case $p\geq 4$. When $p = 3$, we have $$\cot\mu = \tan \left(\pi - \frac \pi 3\right)\frac{\cos\frac \pi 4 + \cos \frac\pi 3 \cos \frac\pi 3}{\sin \frac\pi 3 \sin \frac\pi 3} \geq \tan\left(\pi - \frac {2\pi}{3}\right)\frac{\cos \frac\pi 4 +\cos\frac \pi 3 \cos \frac\pi 3}{\sin \frac \pi 3  \sin \frac\pi 3},$$ which yields the same lower bounds for $\theta$. As before, $\phi \leq \frac \pi 2$, whence the acute quadrilateral as in Lemma \ref{lem:geod_dont_cross}. 
\end{proof}

\begin{lem}\label{lem:zigzag_geodesic_dont_cross}
	Suppose that $q\geq 4$ and that $\gamma_m = 0$ is a point of type $B$ as described above and let $$\mathcal{G}_{j, a\leq k \leq b}^\pm:=\{\lambda_{j,k}^\pm:a\leq k \leq b\}.$$Then geodesics in $\mathcal{G}_{0,0\leq k \leq r-1}^-\cup\mathcal{G}_{q-1,0\leq k \leq r-1}^+$ don't cross geodesics in $\bigcup_{j = 1}^{q-3} \mathcal{G}_{j,0\leq k \leq r-1}^+\cup \bigcup_{j=2}^{q-2}\mathcal{G}_{j,0\leq k \leq r-1}^-$.
\end{lem}

\begin{figure}[h]
	\centering
	\begin{subfigure}[b]{0.45\textwidth}
		\includegraphics[width=\textwidth]{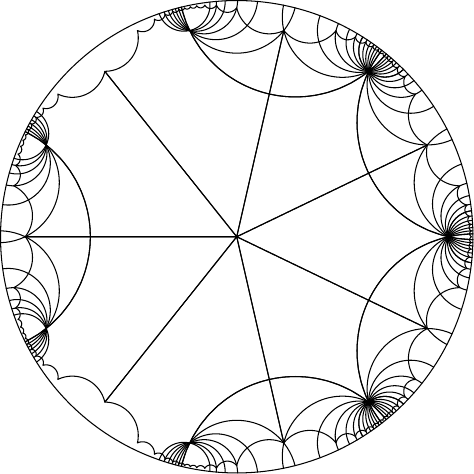}
		\caption{On the left is the group of geodesics $\mathcal{G}_{0,0\leq k \leq r-1}^-\cup\mathcal{G}_{q-1,0\leq k \leq r-1}^+$; on the right is $\bigcup_{j = 1}^{q-3} \mathcal{G}_{j,0\leq k \leq r-1}^+\cup \bigcup_{j=2}^{q-2}\mathcal{G}_{j,0\leq k \leq r-1}^-$.}
		\label{fig:dist_between_geod1}
	\end{subfigure}
	\hfill
	\begin{subfigure}[b]{0.45\textwidth}
		\includegraphics[width=\textwidth]{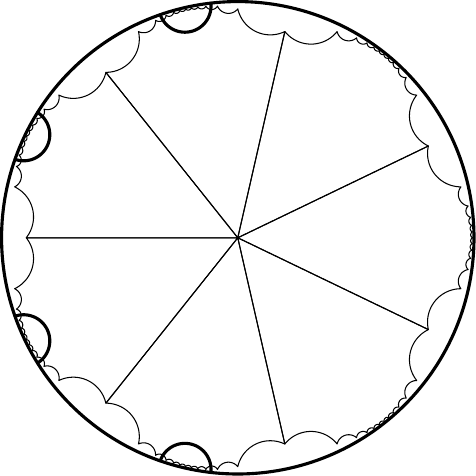}
		\caption{In bold, the four "limiting" geodesics $\lambda_{0,r-1}^-$, $\lambda_{q-1,r-1}^+$, $\lambda_{q-2,r-1}^-$ and $\lambda_{1,r-1}^+$.}
		\label{fig:dist_between_geod2}
	\end{subfigure}
	\caption{}
	\label{fig:dist_between_geod}
\end{figure}

\begin{proof}
	By Lemmas \ref{lem:geod_dont_cross} and \ref{lem:geod_dont_cross2}, there can be no intersection between the geodesics in $\{\lambda_{0,r-1}^-,\lambda_{q-1,r-1}^+\}$ and the ones in $\{\lambda_{1,r-1}^+,\lambda_{q-2,r-1}^-\}$. Since these these geodesics are the outermost in their corresponding sets (see Figure \ref{fig:dist_between_geod}), the geodesics in $\mathcal{G}_{0,0\leq k \leq r-1}^-\cup\mathcal{G}_{q-1,0\leq k \leq r-1}^+$ are disjoint from geodesics in $\bigcup_{j = 1}^{q-3} \mathcal{G}_{j,0\leq k \leq r-1}^+\cup \bigcup_{j=2}^{q-2}\mathcal{G}_{j,0\leq k \leq r-1}^-$.
\end{proof}

In the case where $\gamma_m$ is of type $A$, we have

\begin{lem}\label{lem:zigzag_geodesic_dont_cross_A=0}
	Suppose that $p\geq 4$ and that $\gamma_k = 0$ is a point of type $A$ as described above. Geodesics in $\mathcal{G}_{0,0\leq k \leq r-1}^-\cup\mathcal{G}_{q-1,0\leq k \leq r-1}^+$ don't cross geodesics in $\bigcup_{j = 1}^{q-3} \mathcal{G}_{j,0\leq k \leq r-1}^+\cup \bigcup_{j=2}^{q-2}\mathcal{G}_{j,0\leq k \leq r-1}^-$.
\end{lem}
\begin{proof}
	The proof is the same as in Lemma \ref{lem:zigzag_geodesic_dont_cross}
\end{proof}

Now, suppose that $P$ and $Q$ are face polygons adjacent in an edge $e = [x,y]$, where $e$ has the orientation induced by $P$. Label edges of $P$ counterclockwise from $0$ to $2r-1$ starting with $[x,y]$ and edges of $Q$ counterclockwise from $0$ to $2r-1$ starting with $[x,y]$. Give to $P$ and $Q$ the anticlockwise orientation and label the vertex at the positive end of the edge $e$ by $e^+$. For $ 1 \leq k \leq r-1$, define geodesics $\delta_k^\pm$ and $\varepsilon_k^\pm$ as follows.

\begin{itemize}
	\item If $k$ is odd, $\delta_k^-$ (resp. $\varepsilon_k^-$) is the geodesic going through antipodal vertices of $P$ (resp. $Q$) that goes through the vertex $\left(\frac{k-1}{2}\right)^+$ of $P$ (resp. $Q$) and $\delta_k^+$ (resp. $\varepsilon_k^+$), the one that goes through the vertex $\left(2r-1 - \frac{k-1}{2}\right)^+$ of $P$ (resp. $Q$).\\
	\item If $k$ is even, $\delta_k^-$ (resp. $\varepsilon_k^-$) is the geodesic given by Lemma \ref{lem: infinite_edge_seq} going through opposite edges of $P$ (resp. $Q$) that crosses the edge $\frac{k}{2}$ of $P$ (resp. $Q$) and $\delta_k^+$ (resp. $\varepsilon_k^+$), the one that crosses the edge $2r-1 - \frac{k-2}{2}$ of $P$ (resp. of $Q$).
\end{itemize}

We refer the reader to Figure 
\begin{lem}
	In the configuration described above, the geodesics in $\mathcal{D}^+ : =\{\delta_k^+:1\leq k\leq r-1\}$ don't cross the geodesics in $\mathcal{E}^-:= \{\varepsilon_k^+:1\leq k \leq r-1\}$ and the geodesics in $\mathcal{D}^- :=\{\delta_k^-:1\leq k\leq r-1\}$ don't cross the geodesics in $\mathcal{E}^- :=\{\varepsilon_k^-:1\leq k \leq r-1\}$.
\end{lem}

\begin{proof}
	As in Lemmas \ref{lem:geod_dont_cross} and \ref{lem:geod_dont_cross2}, geodesics $\delta_1^+$ and $\delta_{r-1}^+$ are the outermost geodesics of $\mathcal{D}^+$ and $\varepsilon_1^+$ and $\varepsilon_{r-1}^+$ are the outermost geodesics of $\mathcal{E}^+$. Hence we only need to show that the first pair of geodesics doesn't cross the second pair. Let $d_1^-$ and $d_1^+$ be the negative and positive endpoints of $\delta_1^+$ and $e_{r-1}^-$ and $e_{r-1}^+$ the negative and positive endpoints of the geodesic $\varepsilon$ that goes through vertices $(r-2)^+$ and $(2r-2)^+$ of $Q$. Label vertices in the neighbourhood of $x$ from $y_1$ to $y_{\deg(x)}$ clockwise such that $y = y_1$. Let $\alpha$ be the complete geodesic that goes through $y_3$ and $x$. Since $\alpha$ and $\delta_1^+$ cross, their endpoints are intertwined at infinity. We let $\xi^-$ and $\xi^+$ be the endpoints of $\alpha$ in a way that the cyclic inequality $d_1^- \leq \xi^+ <  d_1^+ \leq \xi^- < d_1^-$ is satisfied. We must be careful that if $\deg (x) = 3$, then $\xi^+ = d_1^-$ and $\xi^- = d_1^+$ because $\alpha = \delta_1^+$, but in all other case, the inequality is strict. Consider the polygon $Q'$ that has $[x, y_2]$ and $[x,y_3]$ as edges and the geodesic $\varepsilon'$ that goes through $y_2$ and the antipodal vertex of $y_2$ in $Q'$. If we let $\eta^-$ and $\eta^+$ be the negative and positive endpoints of $\varepsilon'$, then by Lemma \ref{lem:geod_dont_cross} or \ref{lem:geod_dont_cross2} whether $x$ is of type $A$ or $B$, $\alpha$ and $\varepsilon'$ don't cross. This gives the cyclic inequality $\xi^+ < \eta^+ < \eta^- < \xi^- < \xi^+$. Now, since the geodesics $\varepsilon$ and $\varepsilon'$ both go through $y_2$, we must have that 
	
	As in Lemmas \ref{lem:geod_dont_cross} and \ref{lem:geod_dont_cross2}, geodesics $\delta_1^+$ and $\delta_{r-1}^+$ are the outermost geodesics of $\mathcal{D}^+$ and $\varepsilon_1^+$ and $\varepsilon_{r-1}^+$ are the outermost geodesics of $\mathcal{E}^+$. Hence we only need to show that the first pair of geodesics doesn't cross the second pair. Let $d_1^-$ and $d_1^+$ be the negative and positive endpoints of $\delta_1^+$ and $e_{r-1}^-$ and $e_{r-1}^+$ the negative and positive endpoints of the geodesic $\varepsilon_{r-1}^+$. Our goal is to establish the cyclic inequality $d_1^- < e_{r-1}^+ < e_{r-1}^- < d_1^+ < d_1^-$.
	
	Label vertices in the neighbourhood of $x$ from $y_1$ to $y_{\deg(x)}$ clockwise so that $y = y_1$. 
\end{proof}
Before we show how the previous Lemmas imply lower bounds on the length of geodesics, we introduce the notion of a \emph{contributing polygon} to a $S^f$-admissible path and then attach a geodesic to every contributing polygon. The previous Lemmas will imply that the distance between to such geodesic is positive and this will imply the desired lower bounds.

\begin{figure}[h]
	\centering
	\begin{subfigure}[b]{0.45\textwidth}
		\includegraphics[width=\textwidth]{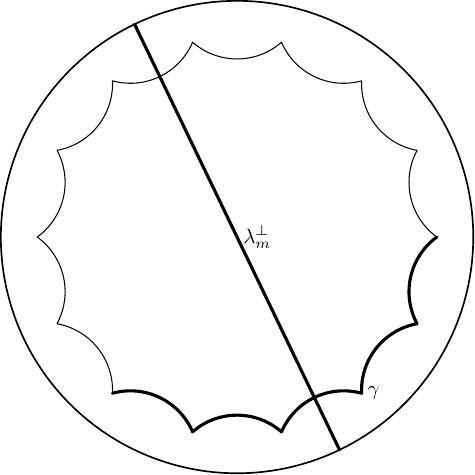}
		\caption{The path $\gamma$ goes through an odd number of sides of $P_m$.}
		\label{fig:lambdaperpvertex1}
	\end{subfigure}
	\hfill
	\begin{subfigure}[b]{0.45\textwidth}
		\includegraphics[width=\textwidth]{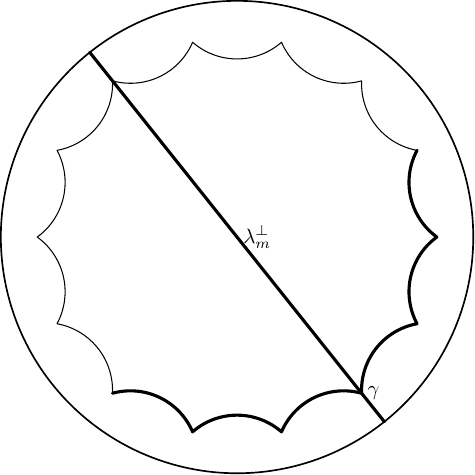}
		\caption{The path $\gamma$ goes through an odd number of sides of $P_m$.}
		\label{fig:lambdaperpvertex2}
	\end{subfigure}
	\caption{}
	\label{fig:lambdaperpvertex}
\end{figure}

\begin{defn}
	Let $\gamma$ be a $S^f$-admissible path. Let $P$ be an adjacent polygon to $\gamma$ (see Definition \ref{def:ajd_interm_polygons}). Let $P'$ be the reflection of $P$ in one of the sides that $\gamma$ shares with $P$. We will say that $\gamma$ is a \emph{contributing polygon to $\gamma$} if
	\begin{itemize}
		\item $P$ shares more than two sides with $\gamma$, or;\\
		\item $\gamma$ goes through exactly one side of $P$ and $P'$, the side through which $P$ and $P'$ intersect, and $P$ is on the left of $\gamma$. 
	\end{itemize}
\end{defn}
In the second case, the condition that $P$ lies on the left of $\gamma$ is arbitrary. In the first case, we will say that $P$ is of type $I$ and in the second, of type $II$. 
\begin{defn}\label{def:geod_attached_to_adjacent_poly}
	Let $\gamma = (\gamma_m)$ be a $S^f$-admissible path and $\{P_j\}$ be its contributing polygons. For every $j$, let $e_1,...,e_{n_j}$ be the edges of $\gamma$ along $P_j$.
	\begin{itemize}
		\item If $n_j$ is even, then $M_j:= \overline{e_{\frac{n_j}{2}-1}}\cap \overline{e_{\frac{n_j}{2}}}$ is a vertex. Define $\lambda_j^\perp$ as the geodesic going through $M_j$ and its opposite vertex in $P_j$ (see Figure \ref{fig:lambdaperpvertex1});\\
		
		\item if $n_j$ is odd, then $M_j := e_{\frac{n_j + 1}{2}}$ is an edge. Define $\lambda_j^\perp$ as the geodesic given by Lemma \ref{lem: infinite_edge_seq} going through $M_j$ and its opposite edge in $P_j$ (see Figure \ref{fig:lambdaperpvertex2}). 
	\end{itemize}
	We also define $p_j :=\lambda_j^\perp \cap M_j$. By our definition of contributing polygons, the $p_j$'s form an increasing sequence along $\gamma$.
\end{defn}

\begin{prop}\label{prop:lambda_perp_dont_cross}
	Let $\gamma$ be a $S^f$-admissible path. For all $j\in \ZZ$, $\lambda_j^\perp$ and $\lambda_{j+1}^\perp$ don't cross. Furthermore, if $\gamma$ is associated with a geodesic $g$, then $\lambda_j^\perp$ and $g$ intersect transversally.
\end{prop}

\begin{proof}
	Write $(P_j)$ for the sequence of contributing polygons to $\gamma$ and fix $j\in \ZZ$. We will denote by $z$ the last vertex of $\gamma$ that goes through $P_j$. There are four cases to consider. Each case depends on whether $P_j$ and $P_{j+1}$ are of type $I$ or $II$.
	\begin{itemize}
		\item \textbf{Case 1:} $P_j$ is of type $I$ and $P_{j+1}$ is of type $I$. Remark that this is only possible if $\deg z \geq 4$. In the notation of Lemma \ref{lem:zigzag_geodesic_dont_cross} (resp. \ref{lem:zigzag_geodesic_dont_cross_A=0}), this corresponds to geodesics $\lambda_j^\perp = \lambda_0^+$ and $\lambda_{j+1}^\perp = \lambda_k^+$ for some $2\leq k \leq q-2$ (resp $2\leq k \leq p-2$) if $z$ is a vertex of type $B$ (resp. of type $A$). These Lemmas show that the geodesics don't intersect.\\
		\item \textbf{Case 2:} $P_j$ is of type $I$ and $P_{j+1}$ is of type $II$. This case is treated as above with $\lambda_{j}^\perp = \lambda_0^+$ and $\lambda_{j+1}^\perp \in \bigcup_{j = 1}^{q-3} \mathcal{G}_{j,0\leq k \leq r-1}^+\cup \bigcup_{j=2}^{q-2}\mathcal{G}_{j,0\leq k \leq r-1}^-$. We conclude that the geodesic don't cross by using Lemma \ref{lem:zigzag_geodesic_dont_cross} of Lemma \ref{lem:zigzag_geodesic_dont_cross_A=0} whether $z$ is of type $B$ or $A$.\\
		\item \textbf{Case 3:} $P_j$ is of type $II$ and $P_{j+1}$ is of type $I$. This case is obtained by the proper rotation of the preceding case. 
		\item \textbf{Case 4:} $P_j$ is of type $II$ and $P_{j+1}$ is of type $II$. The difference here is that if $[y,z]$ is the last edge of $\gamma$ going through $P_j$, then $[y,z]$ may also be an edge of $P_{j+1}$. If this is not the case, we can conclude as all the previous cases. Hence, we will suppose that $[y,z]$ is an edge of both $P_j$ and $P_{j+1}$. Using previous lemmas shows that these geodesics may not cross.
		
		Now, suppose that $\gamma$ comes from a geodesic $g$. By Proposition \ref{prop:ineq_on_5_paths}, we have $\breve{\gamma}_L\prec\gamma \prec \breve{\gamma}_R$ and $\breve{\gamma}_L\ll g \ll\breve{\gamma}_R$. The arguments we've made in the paragraphs above in fact imply that every geodesic $\lambda_j^\perp $ goes through $\breve{\gamma}_L$ and $\breve{\gamma}_R$. Since $g$ has the same endpoints as $\breve{\gamma}_L$ and $\breve{\gamma}_R$ and the strip bounded by these two paths is convex, the $\lambda_j^\perp$'s and $g$ cannot have the same endpoints, hence must cross transversally. 
	\end{itemize}
\end{proof}


\begin{prop}\label{prop:lowebound_on_geod_between_strips_given_by_lambda_perp}
	Let $g$ be a geodesic and $\gamma = \gamma(g)$. Let $\mathcal{R}_m$ be the closed band enclosing $\lambda_m^\perp$ and $\lambda_{m+1}^\perp$ and $g_m := \mathcal{R}_m\cap g$. Let $d_m :=\dist \left(\lambda_m^\perp,\lambda_{m+1}^\perp\right)$. Then $d_m >0$ and $\length(g_m)\geq d_m > 0$.
\end{prop}
\begin{proof}
	In Proposition \ref{prop:existence_of_cvx_paths}, we have found convex paths $\breve{\gamma}_L$ and $\breve{\gamma}_R$ associated with $\gamma = \gamma(g)$ that satisfy $\breve{\gamma}_L\ll g \ll \breve{\gamma}_R$. By Proposition \ref{prop:lambda_perp_dont_cross}, for every $m$, $\gamma$ and $\lambda_m^\perp$ don't have the same endpoints. Hence, endpoints of $\lambda_m^\perp,\lambda_{m+1}^\perp,\breve{\gamma}_L$ and $\breve{\gamma}_R$ are intertwined in $\del\DD$. Let $\mathcal{C}$ be the closed region bounded by $\breve{\gamma}_L$ and $\breve{\gamma}_R$. Remark that $\mathcal{R}_m\cap \mathcal{C}$ is convex and topologically a rectangle. By Proposition \ref{prop:existence_of_cvx_paths}, $g$ crosses the sides $\lambda_m^\perp$ and $\lambda_{m+1}^\perp$. Since $g$ is geodesic, we have $$\length(g\cap \mathcal{R}_m) = \length (g_m) \geq d_m.$$
	
	Since $\lambda_m$ and $\lambda_{m+1}$ do not cross, the latter is indeed strictly positive.
\end{proof}

\begin{prop}\label{prop:closed_geod_give_periodic_paths}
	Let $g$ be a closed geodesic in $X$, $\tilde{g}$ be a lift of $g$ in $\DD$ and $\overline{g}$ be the complete geodesic in $\DD$ defined by $\tilde{g}$. Then $\gamma = \gamma(\overline{g})$ is periodic. Conversely, if $\gamma$ is a periodic $S^f$-admissible path, then $\gamma$ defines a closed geodesic in $X$.
\end{prop}
\begin{proof}
	To prove the first part of the statement, write $\gamma = \gamma(\overline{g}) = (...,P^{-2},P^{-1},P^0,P^1,P^2,...)$. The lift $\tilde{g}$ defines a hyperbolic translation $T_{\tilde{g}}\in G_{p,q,r}$. Let $k$ be such that $P^k = T_{\tilde{g}}(P^0)$ (here, we point out that $P^0$ and $P^k$ must be of the same type). We must show that $P^{k+1} = T_{\tilde{g}}(P^1)$. For this, first remark that $I_{P^k\to P^{k+1}} = I_{T_{\tilde{g}}(P^0)\to T_{\tilde{g}}(P^1)}$. Indeed, we have $T_{\tilde{g}}(I_{P^0\to P^1}) = I_{T_{\tilde{g}}(P^0)\to T_{\tilde{g}}(P^1)} = I_{T_{P^k\to T_{\tilde{g}}(P^1)}} $. Let $\xi^+$ be the positive endpoint of $\overline{g}$. Among the $\deg P^k$ possibilities for this interval, $I_{P^k\to P^{k+1}}$ is the only one that contains $\xi^+$. Hence, to show the desired equality, it suffices to show that $\xi^+\in T_{\tilde{g}}(I_{P^0\to P^1})$, but this is clearly the case since $\xi^+\in I_{P^0\to P^1}$ implies that $\xi^+ = T_{\tilde{g}}(\xi^+)\in T_{\tilde{g}}(I_{P^0\to P^1})$. Since the interval $I_{P^k\to Q}$ uniquely determines the point $Q$, we must have $P^{k+1} = T_{\tilde{g}}(P^1)$, so that $\gamma$ is indeed periodic. 
	
	Conversely, if $\gamma$ is any $S^f$-admissible path, then its code $w = w(\gamma)$ is periodic. We will write $w = (w')^\infty$, where $w' = a^{e_1}b^{f_1}...a^{e_k}b^{f_k}$ is minimal. Then $w'$ defines the hyperbolic transformation $a^{e_1}b^{f_1}...a^{e_k}b^{f_k}$, where here, $a$ is the anticlockwise rotation of $\frac \pi p$ around $A_0$ and $b$ is the clockwise rotation of $\frac{\pi}{q}$ around $B_0$ (we hope that the reader will not be confused by this slight abuse of notation). This defines a closed geodesic in $X$.
\end{proof}

\begin{defn}
Let $\mathcal{C}_\gamma := \{P_j\}$ be the set of contributing polygons to $\gamma = \gamma(g)$, where $g$ is a closed geodesic. The group $\langle \gamma\rangle$ acts on $\mathcal{C}_\gamma$. The \emph{combinatorial length} of $\gamma$ is $$L(\gamma):= \#\{\orb_P:P\in \mathcal{C}_\gamma\},$$
	where $\orb_P$ is the orbit of the polygon $P\in\mathcal{C}_\gamma$.
\end{defn}

\begin{thm}\label{thm:bound_on_length}
	There exists a constant $c = c_{p,q,r}>0$ which depends only on the group $G_{p,q,r}$ for which we have the inequality $$\length(g) \geq cL(\gamma(g))$$ for any closed geodesic $g$.
\end{thm}

\begin{proof}
	First, we point out that the proof of Proposition \ref{prop:lambda_perp_dont_cross} consists in examining all possible confirgurations of contributing polygons (which there are a finite number of) and showing that in each case, $d_m > 0 $ (or equivalently that the geodesics are disjoints in $\overline{\DD}$). Since, for a fixed hyperbolic triplet, $(p,q,r)$, there is only a finite number of cases, we can define $c := \min{d_m}$, where the minimum is taken over all possible configurations.
	
	Now, lift the geodesic $g$ to a complete connected geodesic $\overline{g}\in\DD$. By Proposition \ref{prop:closed_geod_give_periodic_paths}, the path $\gamma = \gamma(\overline{g}) = (...,P^{-2},P^{-1}, P^0, P^1, P^2, ...)$ is periodic. Let $\{P_j\}$ be the set of contributing polygons to $\gamma$ and $(\lambda_j^\perp)$ be the geodesic attached to $P_j$. Through the proof of Proposition \ref{prop:closed_geod_give_periodic_paths}, we've shown that there exists $k\geq 0$ such that $T_{\tilde{g}}(P^n) = P^{n+k}$ for all $n\in \ZZ$ and for the proper (isometric) lift $\tilde{g}$ of $g$. This implies that there exists $k'$ such that $T_{\tilde{g}}(\lambda_j) = \lambda_{j+k'}$ for all $j\in \ZZ$. Hence, the band $\mathcal{B}$ bounded by the geodesics $\lambda_0$ and $\lambda_{k'}$ gives rise to the space $Y := \mathcal{B}/\sim$, obtained by glueing points on $\lambda_0$ to points of $\lambda_{k'}$ via $T_{\tilde{g}}$. In this space, every contributing polygon is represented once; there are precisely $k'$ of them. Hence, $L(\gamma(g)) = k'$. The space $Y$ inherits the hyperbolic metric from $\DD$ (it is isometric to the space given by the subgroup $\langle g \rangle$ of $\pi_1(X)$). Denote by $[\mathcal{R}_j]$ the image in $Y$ of the band $\mathcal{R}_j$ defined in Proposition \ref{prop:lowebound_on_geod_between_strips_given_by_lambda_perp}. Then 
	
	\begin{align*}
		\length(g) &= \sum_{j = 0} ^{k'-1} \length (g|_{[\mathcal{R}_j]}) \\
		&=\sum_{j = 0} ^{k'-1} \length (g|_{\mathcal{R}_j}) \\
		&\geq \sum_{j = 0} ^{k'-1}d = ck' = cL(\gamma(g)).
	\end{align*}
	
\end{proof}

The combinatorial length can in fact be given a more concrete formulation when we know the code of a periodic $S^f$-admissible path. We first define some words that act as "building blocks" of words on the letters $a^e$ and $b^f$.

\begin{defn}
	Let $(p,q,r)$ be a hyperbolic triplet. For $n\in \NN$, we will refer to the words
	\begin{itemize}
		\item $e_A^-(n) := \underbrace{ba^{p-1}ba^{p-1}b...}_{n-1\text{ symbols}}$ ,\\
		\item $e_A^+(n):= \underbrace{b^{q-1}ab^{q-1}ab^{q-1}...}_{n-1\text{ symbols}}$ ,\\
		\item $e_B^-(n) := \underbrace{a^{p-1}ba^{p-1}ba^{p-1}...}_{n-1\text{ symbols}}$ and\\
		\item $e_B^+(n) := \underbrace{ab^{q-1}ab^{q-1}a...}_{n-1\text{ symbols}}$
	\end{itemize}
	as \emph{uni-polygonal words }. We will say that two uni-polygonal $e_1$ and $e_2$ \emph{glue} if the last letter of $e_1$ is different from the first letter of $e_2$. In such a case, we can concatenate these words into $e_1 e_2$. We define $e_A^\pm(1)$ and $ e_B^\pm(1)$ to both be the empty word.
\end{defn}

The words defined above are seen to border only one face polygon (without restriction on the number $n$, they may wind around it multiple times). The number $n$ in the definition above corresponds to the number of sides of such a polygon a path on $\gr$ goes through to have such a code. Remark that given $n_1,n_2 \in\NN$ and a uni-polygonal word $e_A^+(n_1)$, then there exists a unique uni-polygonal word $e_P^-(n_2)$ that glues to $e_A^+(n_1)$. The same goes if we had started with any other type of uni-polygonal word. This motivates the following definition.

\begin{defn}\label{def:multi-polygonal words}
	Let $n_1,n_2,...,n_k\in\NN_{\geq 2}$. We will refer to the words
	\begin{itemize}
		\item $z_A^-(n_1,n_2,...,n_k):= e_A^-(n_1) e_{P_2}^+(n_2-1) e_{P_3}^-(n_3-1)... e_{P_k}^t(n_k-1)$,\\
		\item $z_A^+(n_1,n_2,...,n_k) := e_A^+(n_1) e_{P_2}^-(n_2-1) e_{P_3}^+(n_3-1) ... e_{P_k}^t(n_k-1)$,\\ 
		\item $z_B^-(n_1,n_2,...,n_k) := e_B^-(n_1) e_{P_2}^+(n_2-1) e_{P_3}^-(n_3-1) ... e_{P_k}^t(n_k-1)$ and\\
		\item $z_B^+(n_1,n_2,...,n_k) := e_B^+(n_1) e_{P_2}^-(n_2-1) e_{P_3}^+(n_3-1) ... e_{P_k}^t(n_k-1)$,
	\end{itemize}
	as \emph{zigzags}, where, in each case, the type of point $P_j$ is determined by the integer $n_{j-1}$ and $t = +$ or $-$ depending on the parity of $k$. We will refer to the number $k$ as the \emph{length} of a zigzag $\zeta$, denoted by $L(\zeta)$.
\end{defn}

\begin{figure}[h]
	\centering
	\begin{subfigure}[b]{0.45\textwidth}
		\includegraphics[width=\textwidth]{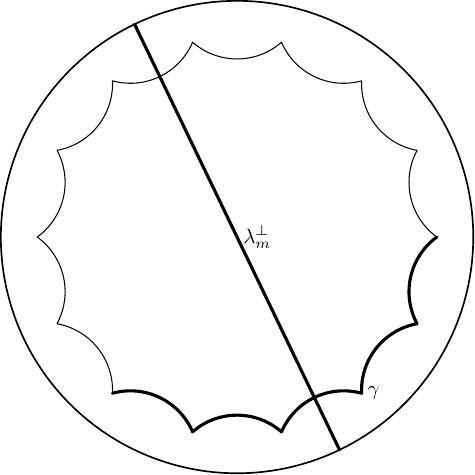}
		\caption{The path $\gamma$ goes through an odd number of sides of $P_m$.}
		\label{fig:example_multipolygonal_words1}
	\end{subfigure}
	\hfill
	\begin{subfigure}[b]{0.45\textwidth}
		\includegraphics[width=\textwidth]{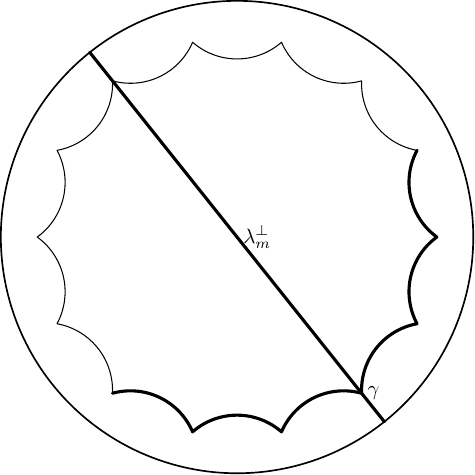}
		\caption{The path $\gamma$ goes through an odd number of sides of $P_m$.}
		\label{fig:example_multipolygonal_words2}
	\end{subfigure}
	\caption{}
	\label{fig:example_multipolygonal_words}
\end{figure}

\begin{exa}
	Say you're given $(n_1,n_2,n_3) = (2,4,3)$ with $(p,q,r) = (3,4,4)$ and you want to form the word $z_A^+(n_1,n_2,n_3)$. Then $e_A^+(n_1) =  b^3$. Since $n_1$ is even, we have $P_2 = B$, so that $e_{P_2}^-(n_2-1) = e_A^-(3) = a^2 b a^2$ and since $n_2$ is even, $P_3 = A$, so that $e_{P_3}^+(n_3 - 1) = e_{A}^+(3) = b^3 a$. Hence, $z_A^+(n_1,n_2,n_3) = z_A^+ (2,4,3) = b^3a^2 b a^2 b^3 a$ (see Figure \ref{fig:example_multipolygonal_words}).
\end{exa}

\begin{rem}\label{rem:length_and_contributing_polygons_for_zigzags}
	Any zigzag $\zeta$ is, by definition a product of uni-polygonal words. To each of those uni-polygonal word, we may attach its associated contributing polygon. From this point of view, the length $L(\zeta)$ counts the number of those contributing polygons.
\end{rem}

If we restrict ourselves to finite words that do not contain the symbols $a^e$ for $2\leq e \leq p-2$ and $b^f$ for $2\leq f\leq q-2$, words we've defined as switches above, then Definition \ref{def:multi-polygonal words} is merely an indexation of these words, i.e. any finite words without switch is represented as a zigzag.

\begin{prop}\label{prop: factorization_in_subwords}
	Let $\gamma$ be a periodic $S^f$-admissible word. Then there exists a finite set of zigzags $\mathcal{Z} = \{\zeta_j\}$ and a finite set of switches $\mathcal{S} = \{\sigma_j\}$ such that $$w(\gamma) = \left(\prod_j\zeta_j\sigma_j\right)^\infty.$$
\end{prop}

\begin{proof}
	We write $\gamma = (...,P^{-2}, P^{-1},P^0,P^1,P^2,...)$ and define $S := \{P^k:P^k \text{ defines a switch.}\}$. By definition, the maximal subwords of $\gamma':= \gamma\setminus S$ are zigzags. Hence, we may select an element $\zeta_0\in \gamma'$. If we order the maximal subwords of $\gamma'$ coming after $\zeta_0$ then those must come periodically due to the periodicity of $\gamma$. Therefore, there exists a minimal $k$ such that $\zeta_0 = \zeta_k$. We define $\mathcal{Z}:=\{\zeta_j:0\leq j\leq k-1\}$. Let $\sigma_j$ be the switch in between $\zeta_j$ and $\zeta_{j+1}$ and define $\mathcal{S}:=\{\sigma_j:0\leq j\leq k-1\}$. It may happen that $S = \varnothing$. In that case, there is only one maximal subword of $\gamma'$, which the whole $\gamma$. In that case, let $k$ be the smallest integer for which $P^{j+k} = P^j$ for all $j\in\ZZ$. Fix a point $P^{j_0}$ which is a vertex of two contributing polygons. Consider the finite sequence $\alpha := (P^{j_0},P^{j_0 + 1},..., P^{j_0 + k })$. Then by construction, $\gamma = \alpha^\infty$ and since $\gamma$ contains no switches, the code $w(\alpha)$ is a zigzag. We can define $\zeta_0 = w(\alpha)$ and we obtain that $w(\gamma) = \zeta_0^\infty$.
	\end{proof}

Remark that in the preceding proof, the set $S'$ may be empty. This happens if and only if $\gamma$ is an infinite product of zigzags, e.g. $(ab)^\infty$. This factorization also relates the length of a periodic path with the length of its switches.

\begin{prop}\label{prop:Length_in_terms_of_zigzags}
	Let $\gamma = \gamma(g)$ be the $S^f$-admissible path associated to the closed geodesic $g$. If the factorization of $w(\gamma)$ given by Proposition \ref{prop: factorization_in_subwords} contains a switch, then $L(\gamma) = \sum L(\zeta_j)$. If $\gamma$ has no switch, then $L(\gamma) = L(\zeta_0) +1$.
\end{prop}

\begin{proof}
	The factorization of $w(\gamma)$ given by Proposition \ref{prop: factorization_in_subwords} can be written as $$ w(\gamma) =\left( e^0_0e^0_1...e^0_{n_0}\sigma_0 e^1_0...\sigma_{k-2}e^{k-1}_0...e^{k-1}_{n_{k-1}}\sigma_{k-1}\right)^\infty,$$
	
	where, for $0\leq j\leq k-1$ and $0\leq i\leq n_j$, $e^j_i$ is a uni-polygonal word in the switch $\zeta_j$. Consider the map $\Phi$ which maps every $e_i^j$ to its associated contributing polygon. Since different uni-polygonal words maps to different contributing polygons, this map is injective.
	
	Let $P_0$ be the contributing polygon associated with $e_0^0$ and $P_1$ be the one associated with $e_{n_{k-1}}^{k-1}$. If $Q$ is a contributing polygon for which $P_0 < Q < P_1$, then the intersection $Q\cap \gamma$ defines a uni-polygonal word which belongs in the factorization of $\gamma$. 
	
	Now, consider the case where $w(\gamma)$ contains a switch, then the first contributing polygon coming after $P_1$ is $T_g(P_0)$, where $T_g$ is the hyperbolic transformation induced by $g$. Indeed, in that case, there is a switch $\sigma_{k-1}$ and the contributing polygons having that switch as a vertex are precisely $P_1$ and $T_g(P_0)$. Hence, we get that $\sum_j L(\zeta_j) =  L(w(\gamma))$. 
	
	If $w(\gamma)$ contains no switch, we claim that there is a contributing polygon coming in between $P_1$ and $T_g(P_0)$. First remark that this is possible if and only if $k = 1$. Since $w(\gamma)$ contains no switch, the composition of uni-polygonal words $e^0_{n_{0}-1}e_0^0$ contains three contributing polygons. Therefore, to compute the length of $w(\gamma)$ in terms of its factorization, we must account for this extra polygon. In particular, $L(w(\gamma)) = L(\zeta_0) +1$.
\end{proof}

We are now ready to describe how one may enumerate conjugacy classes up to a certain length.

Let $\ell_0>0$ represent a geometric length. Say you want to compute all conjugacy classes up to length $\ell_0$ in $G_{p,q,r}$. Define $L_0 := \frac{\ell_0}{c}$ where $c= c_{p,q,r}$ is the constant given by Theorem \ref{thm:bound_on_length}. Enumerate all words in the corresponding alphabet satisfying the inequalities given by Proposition \ref{prop: admissible_iff_some_inequalities} and the concrete limiting words in Table \ref{thm:bound_on_length_first_time} up to combinatorial length $L_0$. Then for every other word $w = w(g)$ for which $L(w)>L_0$, by Theorem \ref{thm:bound_on_length}, we have 

\begin{align*}
	\length{g}&\geq c L(w(g))\\
	&\geq cL_0 = \ell_0.
\end{align*}

Hence all conjugacy classes of geometric length smaller than $\ell_0$ have been enumerated. Then, one must remove the word $w_L$ if it has been generated and remove all generated words for which the associated geometric length is strictly greater than $\ell_0$.
\bibliographystyle{alpha}
\bibliography{Stopping_criterion}
	

\end{document}